\documentclass[11pt,a4paper]{article}
\usepackage[utf8]{inputenc}
\usepackage{amsmath}
\usepackage{amsfonts}
\usepackage{amssymb}
\usepackage{amsmath}
\usepackage{multicol,color}
\usepackage{float}
\usepackage{graphicx}
\usepackage{amssymb}
\usepackage{theorem}
\usepackage{fancyhdr}
\usepackage{tikz}
\usepackage{enumerate}
\usepackage[margin=1.1in]{geometry}
\usepackage{pgfplots}
\usepgflibrary{shapes.geometric}

\pgfplotsset{my style/.append style={axis x line=middle, axis y line= middle, xlabel={$x$}, ylabel={$y$}, axis equal }}
\usetikzlibrary{through,calc}
\usetikzlibrary{calc,intersections}

\newcommand*{\TlYellow}[1][1]{%
\begin{tikzpicture}[every path/.style={thick,fill=lightgray},scale=#1]
\draw (0,0)    circle (0.07);
\draw (60:0.1) circle (0.07);
\draw[fill=yellow] (0.1,0) circle (0.07);}

\def\disp{\displaystyle}
\def\ve{\varepsilon}
\def\dd{\delta}

\def\lm{\lambda}
\def\O{\Omega}

\def\tilde{\widetilde}

\def\ox{\bar{x}}
\def\oy{\bar{y}}

\def\ov{\bar{v}}

\def\gph{\hbox{}}

\def\tto{\rightrightarrows}

\def\Limsup{\mathop{{\rm Lim}\,{\rm sup}}}

\def\hat{\widehat}
\def\Hat{\widehat}
\def\tilde{\widetilde}

\def\Bar{\overline}

\def\ve{\varepsilon}
\def\B{I\!\!B}
\def\h{\hfill\Box}
\def\R{\mathbb{R}}
\def\N{\mathbb{N}}

\def\gph{\mbox{\rm gph}\,}

\def\dom{\mbox{\rm dom}\,}

\def\dist{\mbox{\rm dist}}

\def\O{\Omega}

\def\ph{\varphi}
\def\emp{\emptyset}

\def\oR{\Bar{\R}}
\def\lm{\lambda}

\def\dd{\delta}

\def\al{\alpha}

\def\N{I\!\!N}
\def\th{\theta}

\newtheorem{theorem}{Theorem}[section]
\newtheorem{lemma}[theorem]{Lemma}

\newtheorem{proposition}[theorem]{Proposition}
\newtheorem{definition}[theorem]{Definition}
\theoremstyle{plain}{\theorembodyfont{\rmfamily}
}
\theoremstyle{plain}{\theorembodyfont{\rmfamily}
}
\theoremstyle{plain}{\theorembodyfont{\rmfamily}
}
\theoremstyle{plain}{\theorembodyfont{\rmfamily}
}

\theoremstyle{plain}{\theorembodyfont{\rmfamily}
}

\begin{document}
\begin{center}
{\bf DISCRETE APPROXIMATIONS AND OPTIMALITY CONDITIONS\\ FOR INTEGRO-DIFFERENTIAL INCLUSIONS}\\[3ex]
ABDERRAHIM BOUACH\footnote{Laboratoire LITAN, \'Ecole sup\'erieure en Sciences et Technologies de l'Informatique et du Num\'erique, RN 75, Amizour 06300, Bejaia, Alg\'erie. (abderrahimbouach@gmail.com).} $\;$
TAHAR HADDAD\footnote{Laboratoire LMEPA, Facult\'e des Sciences Exactes et Informatique, Universit\'e Mohammed Seddik Benyahia, Jijel, B.P. 98, Jijel 18000, Alg\'erie (haddadtr2000@yahoo.fr).}$\;\;$
BORIS S. MORDUKHOVICH\footnote{Department of Mathematics, Wayne State University, Detroit, Michigan 48202, USA (aa1086@wayne.edu). Research of this author was partly supported by the US National Science Foundation under grant DMS-1808978, by the Australian Research Council Discovery Project DP-190100555, and by Project 111 of China under grant D21024.}
\end{center}
\small{\bf Abstract.} This paper addresses a new class of generalized Bolza problems governed by nonconvex integro-differential inclusions with endpoint constraints on trajectories, where the integral terms are given in the general (with time-dependent integrands in the dynamics) Volterra form. We pursue here a threefold goal. First we construct well-posed approximations of continuous-time integro-differential systems by their discrete-time counterparts with showing that any feasible solution to the original system can be strongly approximated in the $W^{1,2}$-norm topology by piecewise-linear extensions of feasible discrete trajectories. This allows us to verify in turn the strong convergence of discrete optimal solutions to a prescribed local minimizer for the original problem. Facing intrinsic nonsmoothness of original integro-differential problem and its discrete approximations, we employ appropriate tools of generalized differentiation in variational analysis to derive necessary optimality conditions for discrete-time problems (which is our second goal) and finally accomplish our third goal to obtain necessary conditions for the original continuous-time problems by passing to the limit from discrete approximations. In this way we establish, in particular, a novel necessary optimality condition of the Volterra type, which is the crucial result for dynamic optimization of integro-differential inclusions.\\[0.5ex]
\vspace*{-0.2in}

\section{Introduction and Problem Formulation}\label{sec:intro}
\setcounter{equation}{0}\vspace*{-0.1in}

In this paper, we focus on the study of a general dynamic optimization problem governed by constrained {\em integro-differential inclusions}. The problem is labeled as $(P)$ and is formulated as follows:
\begin{equation}\label{Bolza-functional}
\text{minimize}\quad J[x]:=\varphi\big(x(T)\big)+\int_{0}^{T}l\big(t,x(t),\dot{x}(t)\big)\,\mathrm{d}t,
\end{equation}
over absolutely continuous trajectories $x: [0,T]\rightarrow \mathbb{R}^{n}$ of the {\em integro-differential inclusion} IDI
\begin{equation}\label{inclusion}
(D_{F,g}) :\;\;\left\{
\begin{array}{l}
\dot{x}(t)\in F\big(t,x(t)\big)+\displaystyle\int_{0}^{t}g\big(t,s,x(s)\big)\,\mathrm{d}s\quad \text{a.e.}\;\;t\in [0,T],\vspace*{0.3cm}\\
x(0)=x_0,
\end{array}\right.
\end{equation}
on the fixed time interval $I:=[0,T]$ subject to the endpoint constraint
\begin{equation}\label{endpoint constraints}
x(T)\in \Omega\subset\mathbb{R}^{n}.
\end{equation}
While the case of {\em Lipschitzian} differential inclusions with $g\equiv 0$ in \eqref{inclusion} has been broadly investigated in the literature from the viewpoint of deriving necessary optimality conditions (see, e.g., the books \cite{clsw,Boris,vinter} and the references therein), the authors are not familiar with any publication concerning integro-differential inclusions of type \eqref{inclusion}. Note that we consider here the most challenging case of the {\em general Volterra form} in \eqref{inclusion}, where the integrand $g$ depends on the time variable $t$.

Let us now list the main assumptions used in the paper.
\begin{enumerate}
\item[$ (\mathcal{H}^{F}) $] $ F: I\times \mathbb{R}^{n}\rightrightarrows\mathbb{R}^{n} $ is a set-valued mapping/multifunction with nonempty closed values, Hausdorff continuous for a.e.
$ t\in I $ uniformly
in $ x\in\mathbb{R}^{n} $, and satisfying the following conditions:
\begin{enumerate}
\item[$ (\mathcal{H}^{F}_{1}) $] There exists nonnegative constant $m_{F}$ for which we have
\begin{equation*}
F(t,x)\subset m_{F}\B\quad\text{whenever}\;\;t\in I,\;\;x\in\mathbb{R}^{n},
\end{equation*}
where $ \B $ denotes the closed ball of $ \mathbb{R}^{n} $ centered at the origin with radius $r=1$.
\item[$ (\mathcal{H}^{F}_{2}) $] There exists nonnegative constant $ l_{F} $ ensuring the estimate
\begin{equation}\label{lipshitz:F}
F(t,x_{1})\subset F(t,x_{2})+l_{F}\Vert x_{1}-x_{2} \Vert\B\quad\text{for all}\;\;x_{1}, x_{2}\in\mathbb{R}^{n},\;t\in I.
\end{equation}
\end{enumerate}
\item[$ (\mathcal{H}^{g}) $] $ g : I^2\times\mathbb{R}^{n} \rightarrow \mathbb{R}^{n} $ is measurable in $(s,t)\in I^2$ for each $x\in\mathbb{R}^{n}$) and such that
\begin{enumerate}
\item[$ (\mathcal{H}^{g}_{1}) $] There exists a nonnegative constant $ \beta $ with
\begin{equation*}
\lVert g(t,s,x)\rVert\leq \beta(1+\lVert x \rVert)\quad\text{for all}\;\;(t,s)\in Q_{\Delta}\;\;\text{and }\;\;x\in\mathbb{R}^{n}.
\end{equation*}
\item[$(\mathcal{H}^{g}_{2}) $] The mapping $ x\mapsto g(t,s,x) $ is continuously differentiable for all $ (t,s)\in I^{2} $, and for each real $ \eta>0 $ there exists a nonnegative constant $ \alpha $
such that
\begin{equation*}
\lVert \nabla\,g(t,s,x) \rVert\leq \alpha\quad\text{for all}\;\;(t,s)\in Q_{\Delta}\;\;\text{and}\;\;x\in\eta\B,
\end{equation*}
\end{enumerate}
where $Q_{\Delta}$ stands for the set
$$
Q_{\Delta}:=\big\{(t,s)\in I^{2}\;\big|\;s \leq t\big\}.
$$
\item[$ (\mathcal{H}^{\varphi,l}) $] The terminal cost $ \varphi : \mathbb{R}^{n}\to\oR$ is lower semicontinuous (l.s.c.), while and the
running cost $ l : [0,T]\times\mathbb{R}^{2n}\to\oR $ is l.s.c.\ with respect to all but time variables being continuous with respect to $t$ and being majorized by a summable function on $[0,T]$ along the reference curves. Furthermore, we assume that $l(t,\cdot)$ is bounded from below on bounded sets for a.e.\ $t\in[0,T]$.
\item[$ (\mathcal{H}^{\Omega}) $] The set $ \Omega $ is closed and bounded in $ \mathbb{R}^{n} $.
\end{enumerate}
Observe that condition \eqref{lipshitz:F} is equivalent to the {\em uniform Lipschitz continuity}
\begin{equation*}
\mathrm{haus}\,\big(F(t,x_{1}),F(t,x_{2})\big)\leq l_{F}\Vert x_{1}-x_{2} \Vert\;\mbox{ for all }\;x_{1},x_{2}\in\mathbb{R}^{n},
\end{equation*}
of $ F(t,\cdot) $ with respect to the {\em Pompeiu-Hausdorff metric} $ \mathrm{haus}\,(\cdot,\cdot) $ on the space of nonempty and compact subsets of $ \mathbb{R}^{n}$, whenever $t\in[0,T]$.

To handle efficiently the Hausdorff continuity of $ F(\cdot,x) $ for a.e. $ t\in [0,T] $, define the {\em averaged modulus of continuity} for the multifunction
 $ F $ in $ t\in [0,T] $ while $ x\in\mathbb{R}^{n} $ by
\begin{equation}\label{tau}
\tau(F,h):=\int_{0}^{T}\sigma(F,t,h)\,\mathrm{d}t,
\end{equation}
where $ \sigma(F,t,h):=\sup\big\{\omega(F,x,t,h)\;\big|\;x\in\mathbb{R}^{n}\big\}$ with
$$
\omega(F,x,t,h):=\sup\Big\{\mathrm{haus}\,\big(F(t_{1},x),F(t_{2},x)\big)\;\Big|\;t_{1},t_{2}\in\Big[t-\frac{h}{2},t+\frac{h}{2}\Big]\cap[0,T]\Big\}.
$$
It is proved in \cite{Dontchev} that if $F(\cdot,x)$ is Hausdorff continuous for a.e. $ t\in [0,T] $ uniformly in $ x\in\mathbb{R}^{n} $, then we have $ \tau(F,h)\rightarrow 0 $ as
$h\rightarrow 0$.

Our approach to obtain necessary optimality conditions for problem $(P)$ is based on the {\em method of discrete approximations}, which is important for its own sake, not just for deriving optimality conditions. This method was first developed in \cite{m95} for dynamic optimization problems governed by Lipschitzian differential inclusions with finite-dimensional state spaces and then extended in \cite{Boris} to evolution inclusions in infinite dimensions, systems with delays, and of neutral type; see also the references therein. Besides establishing necessary optimality conditions in continuous-time systems, this method has its values in constructing efficient approximations of feasible and optimal solutions to the original systems by their discrete counterparts and deriving necessary conditions for the latter problems, which can be viewed as {\em suboptimality} (almost optimality) conditions for continuous-time problems of dynamic optimization. Let us emphasize that the {\em Lipschitz continuity} of the velocity mapping $F$ with respect to the state variables plays an essential role in implementing the method of discrete approximations, as well as in other known methods to derive necessary optimality conditions for differential inclusions employed in the aforementioned publications.

More recently, the method of discrete approximations has been developed for problems of dynamic optimization governed by {\em highly non-Lipschitzian} (even discontinuous) differential inclusions associated with (Moreau) {\em sweeping processes} of the type
\begin{equation}\label{sw}
-\dot{x}(t)\in N_{C(t)}\big(x(t)\big)\quad \text{a.e.}\;\;t\in [0,T],\quad x(0)=x_0\in C(0),
\end{equation}
where $C(t)$ is a moving set parameterized by some control actions, with controls that can also enter the perturbed dynamics. Among many publications in this direction for convex and nonconvex sets $C(t)$, we mention the representative papers \cite{cm} and \cite{cmnn} with the further references therein. In our paper \cite{bhm}, the method of discrete approximations is applied to deriving necessary optimality conditions for optimal control problems with an integro-differential dynamics, where the controlled Volterra term added to \eqref{sw} does not depend on time. Due to the non-Lipschitzian nature of the controlled sweeping processes, the machinery of discrete approximations is quite different from the Lipschitzian case and relies on {\em second-order variational analysis}; namely, on the computation of the second-order subdifferentials of mappings associated with the sweeping dynamics. The reader is referred to \cite{bhm,cm,cmnn} and the recent book \cite{m24} for more details.

The rest of the paper is organized as follows. Section~\ref{sec:2va} contains the needed preliminaries from generalized differentiation in variational analysis. The next Section~\ref{idi} presents useful results, some of which are new, revolving around well-posedness of IDIs under the imposed assumptions.

In Section~\ref{sec:discrete}, we show that any feasible solution of the integro-differential inclusion \eqref{inclusion} can be constructively approximated by feasible solutions of the discretized inclusions such that their piecewise-linear extensions on the entire interval $[0,T]$ strongly converge to the given solution of \eqref{inclusion} in the $W^{1,2}$-norm topology on $[0,T]$. This result is applied below to deriving necessary optimality conditions for $(P)$ while it has an independent theoretical and numerical interest.

Section~\ref{sec:loc-min} deals with a new type of local minimizers for the integro-differential problem $(P)$ and its relaxation based on the fundamental features of integro-differential systems. Using these continuous-time properties and the results of Section~\ref{sec:discrete}, we establish in Section~\ref{sec:disc-opt} the well-posedness of discrete approximation problems   $(P_k)$, $k\in\N$, for the nonconvex dynamic optimization problem $(P)$ under consideration and strong convergence of discrete optimal solutions of discrete problems to the designated local minimizer of constrained integro-differential inclusions.

Section~\ref{sec:disc-nec} is devoted to deriving necessary conditions for optimal solutions to the discrete-time problems $(P_k)$ whenever $k\in\N$. Note that the discrete problems $(P_k)$ are intrinsically nonconvex even for smooth and convex initial data due to the (nonconvex as a rule) graphical constraints coming from the discretized integro-differential inclusion in \eqref{inclusion}. To handle such constraints in discrete-time and then continuous-time systems, we employ the robust generalized differential constructions initiated by the third author and reviewed in Section~\ref{sec:2va}.

The subsequent Section~\ref{sec:bolza-nec} is the culmination of the paper. Here we establish novel necessary optimality conditions for designated local minimizers of the original problem  $(P)$ governed by constrained integro-differential inclusions by passing to the limit from the obtained conditions for discrete approximations by using the convergence of their optimal optimal and the imposed Lipschitz continuity of the velocity ma[ping $F$. The concluding Section~\ref{sec:concl} summarizes the main achievements of the paper and discusses some directions of the future research. $\h$\vspace*{-0.2in}

\section{Preliminaries from Variational Analysis}\label{sec:2va}\vspace*{-0.1in}

In this section, we first overview some fundamental constructions of generalized differentiation for sets, set-valued mappings, and extended-real-valued functions employed in what follows. Then we present some preliminary results on discrete and integro-differential systems broadly used below.

To proceed, for any set $\Omega\subset\R^n$ locally closed around $\ox\in\Omega$, consider the (Euclidean)
distance function dist$(x;\Omega):=\inf_{y\in\Omega}\|x-y\|$ and define the projection operator
$\Pi\colon\R^n\tto\Omega$ by
\begin{equation}\label{proj}
\Pi(x;\Omega):=\big\{w\in\R^n\;\big|\;\|x-w\|=\dist(x;\Omega)\big\},\quad x\in\R^n.
\end{equation}
The {\em proximal normal cone} to $\Omega$ at $\ox$ is given by
\begin{equation*}
N^P_{\Omega}(\ox):=\big\{v\in\R^n\;\big|\;\mbox{ such that }\;\ox\in\Pi_\Omega(\ox+\alpha v)\big\}\;\mbox{ if }\;\ox\in\O,
\end{equation*}
and $\Pi_\Omega(\ox):=\emptyset$ otherwise: see
\cite{clsw,Mord,rw} for equivalent descriptions and further references. It is easy to see that $N^P_{\Omega}(\ox)$ is a convex cone, which may not be generally closed. Given now a nonempty set $\Omega\subset\R^n$ locally closed around $\ox$, consider the associated projection operator \eqref{proj} and define, as in the original paper \cite{m76}, the (basic/limiting/Mordukhovich) {\em normal cone} to $\Omega$ at $\ox$ by
\begin{equation}\label{nor-cone}
N_\Omega(\ox):=\big\{v\in\R^n\;\big|\;\exists\,x_k\to\ox,\;w_k\in\Pi_\O(x_k),\;\al_k\ge 0\;\mbox{ with }\;\al_k(x_k-w_k)\to v\big\},
\end{equation}
for $\ox\in\Omega$, and $N_\O(\ox):=\emptyset$ if $\ox\notin\Omega$. Note that the normal cone \eqref{nor-cone} is {\em nontrivial} $N_\O(\ox)\ne\{0\}$  at any boundary point of $\O$,  while it may be nonconvex even for simple nonconvex sets as, e.g., for $\O:=\{(x,y)\in\R^n\;|\;y\ge-|x|\}$ and $\O:=\{(x,y)\in\R^n\;|\;y=|x|\}$. Despite the nonconvexity, the basic normal cone the associated coderivative and subdifferential constructions for set-valued mappings and extended-real-valued functions exhibit {\em full calculi}, which are richer than for convex-valued counterparts. This is due to the {\em variational/extremal principles} of variational analysis; see the books \cite{Mord,b3,rw} and their bibliographies. The following useful result is taken from
\cite[Theorem~1.97]{Mord}:
\begin{equation}\label{cone:dist}
N_\O(\ox)={\rm cone}\big[\partial\dist(\ox;\O)\big]\;\mbox{ for all }\;\ox\in\O.
\end{equation}
where `cone' stands for the conic hull of the set in question.

Next we consider a multifunction $F\colon\R^n\tto\R^m$ with the domain and graph given by
\begin{equation*}
\dom F:=\big\{x\in\R^n\;\big|\;F(x)\ne\emp\big\}\;\mbox{ and }\;\gph F:=\big\{(x,y)\in\R^n\times\R^m\;\big|\;y\in F(x)\big\},
\end{equation*}
respectively. The {\em coderivative} $D^{\ast}F(\bar{x},\bar{y})\colon\mathbb{R}^m\rightrightarrows \mathbb{R}^{n} $ of $F$ at $ (\bar{x},\bar{y})$ is defined by
\begin{equation}\label{cod}
D^{\ast}F(\bar{x},\bar{y})(u):=\big\{\upsilon\in\mathbb{R}^{n}\;\big|\;(\upsilon,-u)\in N_{{\rm\small gph}\,F}(\bar{x},\bar{y})\big\}\;\mbox{ for all }\; u\in\mathbb{R}^m,
\end{equation}
where $\oy=F(\ox)$ is skipped if $F\colon\R^n\to\R^m$ is single-valued. If in the latter case $F$ is continuously differentiable around $ \bar{x} $, then we have the representation
\begin{equation*}
D^{\ast}F(\bar{x})(u):=\big\{\nabla F(\bar{x})^{\ast}u\big\}\;\mbox{ for all }\; u\in\mathbb{R}^m,
\end{equation*}
where the $A^*$ indicates the adjoint/transpose matrix of $A$.  It is important to observe that if the basic normal cone \eqref{nor-cone} is replaced in \eqref{cod} by the (Clarke) {\em convexified normal cone} \cite{clsw} represented for any (closed) set $\Omega$ as the convex closure of \eqref{nor-cone} by
\begin{equation}\label{clarke}
\Bar{N}_\O(\ox):={\rm clco}\,N_\O(\ox),\quad\ox\in\O,
\end{equation}
where `co' stands for the convex hull while `cl' signifies the closure of the given set, then the corresponding coderivative may be extremely large; in particular, when the set  is a graph of some Lipschitzian multifunction. For instance, in the case where $\O=\gph|x|$ and $\ox=(0,0)\in\R^2$, we have that $\Bar{N}_\O(\ox)=\R^2$, and this happens not just in this simple example. It follows from the results by Rockafellar \cite{r85} (see also \cite[Theorem~1.46]{Mord}), that the convexified normal cone \eqref{clarke} of any multifunction $F\colon\R^n\tto\R^m$ whose graph is locally homeomorphic around $(\ox,\oy)\in\gph F$ to the graph of a locally Lipschitzian mapping $f\colon\R^n\times\R^m\to\R^d$ around the corresponding point, is actually a {\em linear subspace}  of dimension $d\ge m$ in $\R^n\times\R^m$, where $d=m$ {\em if and only if} $f$ is {\em smooth} around the point in question. Furthermore, the class of such {\em graphically Lipschitzian manifolds} is fairly broad including---besides graphs of single-valued locally Lipschitzian mapping---also graphs of maximal monotone operators, subgradient mappings for extended-real-valued convex, concave, saddle, and prox-regular mappings; see \cite{Mord,r85,rw} for more details and discussions. This shows that using the convexification in \eqref{clarke} is not appropriate for graphically Lipschitzian multifunctions, which unavoidably emerge in the study of differential and integro-differential inclusions as in \eqref{inclusion}.

Given further an extended-real-valued function $ \varphi:\mathbb{R}^{n}\to \bar{\mathbb{R}}:=(-\infty,\infty]$ finite at $\ox$, the {\em subdifferential} of $ \varphi $ at $ \bar{x} $ is defined geometrically by
\begin{equation}\label{sub}
\partial\varphi(\bar{x}): = \big\{\upsilon\in\mathbb{R}^{n}\;\big|\; (\upsilon,-1)\in N_{{\rm\small epi}\,\ph}\big(\bar{x},\varphi(\bar{x})\big)\big\},
\end{equation}
while admitting various analytic representations given in \cite{Mord,b3,rw}. It has been well recognized in variational analysis (see, e.g.,\cite[Theorem~1.22]{b3}) that the Lipschitz continuity of $\ph\colon\R^n\to\oR$ around $\ox$ with modulus $ l_{\ph}$ ensures that
\begin{equation}\label{e:sub}
\partial\ph(\bar{x})\ne\emp\;\mbox{ and }\;\Vert\upsilon\Vert\leq l_{\ph}\;\mbox{ for all }\;\upsilon\in \partial\ph(\bar{x}).
\end{equation}

In fact,  this result can be deduced from the next theorem taken from \cite[Theorem~5.11]{b5}, which provides complete characterizations of local Lipschitzian behavior of closed-graph multifunctions. The equivalences below are expressed in terms of the coderivative \eqref{cod} and are known as the {\em Mordukhovich criterion}; see \cite[Theorem~9.40]{rw}.\vspace*{-0.1in}

\begin{theorem}[\bf coderivative characterizations of Lipschitzian multifunctions]\label{est:cod}
Let $F\colon\R^n\tto\R^m$ be of closed graph and bounded around $ \bar{x} $ with $F(\bar{x})\ne\emp$. Then each of the following conditions is necessary and sufficient for $F$ to be locally Lipschitzian around this point:\vspace*{-0.05in}

{\bf(i)} There exist a neighborhood $ U $ of $ \bar{x} $ and a constant
 $l\geq 0 $ such that
\begin{equation*}
\sup\,\{\Vert u\Vert\,\in D^{*}F(x,y)(\upsilon)\}\le l\Vert\upsilon\Vert\;\mbox{ for all }\;x\in U,\;y\in F(x),\;\mbox{ and }\;\upsilon\in\mathbb{R}^{m}.
\end{equation*}

{\bf(ii)} We have
$D^{*}F(\bar{x},\bar{y})(0)=\{0\}\;\mbox{ whenever }\;\bar{y}\in F(\bar{x})$.
\end{theorem}\vspace*{-0.2in}

\section{Existence result for IDIs}\label{idi}\vspace*{-0.1in}

The main goal of this section is to show that the imposed assumptions on the initial data of integro-differential inclusion $(D_{F,g})$ ensure the existence of solutions to \eqref{inclusion} with desired estimates. We also present some auxiliary results of their own interest. First we recall a discrete version of the classical Gronwall's inequality taken from \cite{clar}.\vspace*{-0.1in}

\begin{proposition}[discrete Gronwall's inequality]\label{granwal}
Let $ e_{n},\rho_{n},\gamma_{n},\sigma_{n}\ge 0$ be such that
\begin{equation*}
e_{n+1}\leq \sigma_{n}+\rho_{n}\sum\limits_{i=0}^{n-1}e_{i}+(1+\gamma_{n})e_{n}\;\mbox{ for all }\;n\in\N.
\end{equation*}
Then whenever $ n\in\N $ we have the estimate
\begin{equation*}
e_{n}\leq\Big(e_{0}+\sum\limits_{i=0}^{n-1}\sigma_{i}\Big)\exp\Big(\sum\limits_{i=0}^{n-1}(i\rho_{i}+\gamma_{i})\Big) .
\end{equation*}
\end{proposition}

Now we consider several lemmas related to Gronwall's inequalities.\vspace*{-0.1in}

\begin{lemma}[\bf discrete recurrent estimates]\label{granwal1}
Let $ x_{j},a_{j},b_{j},c_{j}\ge 0$ be such that
\begin{equation}\label{g1}
x_{j}\leq c_{j}+b_{j}\sum\limits_{i=j+1}^{k}x_{i+1}+(1+a_{j})x_{j+1}\;\mbox{ for all }\;k\in\N\;\mbox{ and }\;j=0,\ldots,k-1.
\end{equation}
If $ x_{k+1}=0 $, then whenever $ k\in\N$ we have the estimates
\begin{equation*}
x_{j+1}\leq\Big(x_{k}+\sum\limits_{i=j+1}^{k-1}c_{i}\Big)\exp\Big(\sum\limits_{i=j+1}^{k-1}\big((i-j-1)b_{i}+a_{i}\big)\Big)\;\mbox{ for all }\;j=0,\ldots,k-2.
\end{equation*}
\end{lemma}\vspace*{-0.1in}
{\bf Proof}. Fixed any $ k\in\N$ and consider the following sequences:
\begin{equation*}
\begin{aligned}
&u_{k-j}=x_{j}\;\text{ for all }\;j=0,\ldots,k\;\text{ with }\;u_{-1}=x_{k+1},\\
&v_{k-j-1}=a_{j},\;w_{k-j-1}=b_{j},\;\gamma_{k-j-1}=c_{j}\;\text{ for all }\;j=0,\ldots,k-1.
\end{aligned}
\end{equation*}
Then, it follows from \eqref{g1} that
\begin{equation*}
u_{k-j}\leq \gamma_{k-j-1} + w_{k-j-1}\sum\limits_{i=0}^{k-j-1}u_{i-1}+(1+v_{k-j-1})u_{k-j-1}.
\end{equation*}
Putting $ n=k-j-1 $, the above inequality and $ u_{-1}=0 $ ensure that
 \begin{equation*}
u_{n+1}\leq \gamma_{n} + w_{n}\sum\limits_{i=0}^{n-1}u_{i} + (1+v_{n})u_{n}.
\end{equation*}
Applying Proposition~\ref{granwal} brings us to
\begin{equation*}
u_{n}\leq\Big(u_{0}+\sum\limits_{i=0}^{n-1}\gamma_{i}\Big)\exp\Big(\sum\limits_{i=0}^{n-1}\big(iw_{i}+v_{i}\big)\Big),
\end{equation*}
which is equivalent in turn to
\begin{equation*}
x_{j+1}\le\Big(x_{k}+\sum\limits_{i=j+1}^{k-1}c_{i}\Big)\exp\Big(\sum\limits_{i=j+1}^{k-1}\big((i-j-1)b_{i}+a_{i}\big)\Big)\;\mbox{ for all }\;j=0,\ldots,k-2,
\end{equation*}
and therefore justifies the claim of the lemma. $\h$\vspace*{0.1in}

The next lemma obtained in \cite[Lemma~3.2]{bou} can be viewed as an extended continuous-time Gronwall's inequality, which is needed to study integro-differential systems.\vspace*{-0.1in}

\begin{lemma}[\bf extended Gronwall's inequality]\label{granwal:cont}
Given a continuous-time interval $[T_0,T]$, let $ \rho:[T_0,T]\rightarrow\mathbb{R}_{+}$  be an absolutely continuous function, and let $ b_{1}, b_{2}, a :[T_0,T]\rightarrow \mathbb{R}_{+} $ be Lebesgue integrable functions on $[T_0,T]$. Assume that
\begin{equation*}
\dot{\rho}(t)\leq a(t)+b_{1}(t)\rho(t)+b_{2}(t)\int_{T_{0}}^{t}\rho(s)\,\mathrm{d}s \hspace{0.3cm} a.e.\,\, t\in[T_0,T].
\end{equation*}
Then whenever $ t\in[T_0,T]$, we have the estimate
\begin{equation*}
\begin{aligned}
\rho(t)\leq \rho(T_{0})\,\exp\Big(\int_{T_{0}}^{t}\big(b(\tau)+1\big)\,\mathrm{d}\tau\Big) + \int_{T_{0}}^{t}a(s)\,\exp\Big(\int_{s}^{t}\big(b(\tau)+1\big)\,\mathrm{d}\tau\Big)\,\mathrm{d}s,
\end{aligned}
\end{equation*}
where $ b(t):=\max\{b_{1}(t),b_{2}(t)\} $ for all $ t\in[T_0,T]$.
\end{lemma}\vspace*{-0.1in}

Consider further a sequence of partitions $\{\pi_{k}\}_{k\in\N}$, not necessarily uniform, of the continuous-time interval $I=[0,T]$  given by
\begin{equation}\label{mesh}
\Delta_k:=\big\{0=t^k_0<t^k_1<\ldots<t^k_k=T\big\}\;\mbox{ with }\;h^k_j:=t^{k}_{j+1}-t^k_j\;\mbox{ and }\;\max_{j=0,\ldots,k-1}{h^k_j}\le h_k:=\frac{T}{k}
\end{equation}
and recall the following useful result from \cite[Lemma~3.3]{Carstensen}.\vspace*{-0.1in}

\begin{lemma}[\bf linear operators over discrete meshes]\label{mu:lemma}
For a real number $r\in[1,\infty)$ and some partition $\Delta_k$ from \eqref{mesh}, define the sequence of linear operators $ \mu_{k}: L^{r}(I,\mathbb{R}^{n})\rightarrow L^{r}(I,\mathbb{R}^{n}) $ by
\begin{equation*}
\mu_{k}\,y:=\dfrac{1}{t^{k}_{j+1}-t^{k}_{j}}\int\limits_{t^{k}_{j}}^{t^{k}_{j+1}}y(t)\,\mathrm{d}t\;\mbox{ for all }\;k\in\N.
\end{equation*}
Then the operators $\mu_k$ are $L^{1}$uniformly bounded as $k\in\N$, and we have the convergence $\mu_{k}\,y\rightarrow y$ in $ L^{r}(I,\mathbb{R}^{n}) $ as $ k\rightarrow\infty $ whenever $ y\in L^{r}(I,\mathbb{R}^{n})$.
\end{lemma}\vspace*{-0.1in}

Now we are in a position to establish the {\em existence of solutions}  for the integro-differential inclusion $(D_{F,g})$ in \eqref{inclusion} and ensure, in particular, the existence of its solution $x(\cdot)$ belonging to the space $W^{1,2}$ with constructive uniform estimates of the solution and its velocity via the given data of \eqref{inclusion}.\vspace*{-0.1in}

\begin{theorem}[\bf existence of solutions for integro-differential inclusions]\label{exist:solu}
Consider $\mathrm{IDI}$ in $(D_{F,g})$ under the assumptions in $(\mathcal{H}^{F})$ and $(\mathcal{H}^{g})$. Then for any initial point $ x_{0} $, system \eqref{inclusion} admits at least a  trajectory $x(\cdot)\in W^{1,2}([0,T],\mathbb{R}^{n})$. Denoting further
\begin{equation*}
M_{1}:=\Big(1+\Vert x_{0} \Vert +\dfrac{m_{F}}{\beta+1}\Big)\exp\big(T(\beta+1)\big)\quad\text{and}\quad M_{2}:=m_{F}+\beta T\,M_{1},
\end{equation*}
we have the uniform estimates
\begin{equation}\label{bound:solution}
1+\Vert x(t) \Vert\leq M_{1}\;\mbox{ and }\;\Vert \dot{x}(t) \Vert\leq M_{2}\;\text{ on }\;[0,T].
\end{equation}
\end{theorem}\vspace*{-0.1in}
{\bf Proof}. It is easy to check that all the assumptions of \cite[Theorem~3.3]{Gaouir} are satisfied, and thus we deduce the existence statement from the aforementioned result.
It remains to verify the estimates in \eqref{bound:solution}, which play a fundamental role in our subsequent device. To proceed, let $ x(\cdot) $ be the unique solution of $(D_{F,g})$, i.e.,
\begin{equation*}
\begin{aligned}
\dot{x}(t)\in F\big(t,x(t)\big)+\int_{0}^{t}g\big(t,s,x(s)\big)\,\mathrm{d}s\quad\text{a.e.}\;\;t\in [0,T],\;\;x(0)=x_{0}.
\end{aligned}
\end{equation*}
It follows from the imposed assumptions $ (\mathcal{H}^{F}_{1}) $ and $ (\mathcal{H}^{g}_{1}) $ that
\begin{equation}\label{e:rho}
\begin{aligned}
\Vert\dot{x}(t)\Vert\leq m_{F}+\beta\int_{0}^{t}(1+\Vert x(s)\Vert)\,\mathrm{d}s\;\mbox{ for a.e. }\;t\in [0,T].
\end{aligned}
\end{equation}
Setting $\rho(t):=1+\Vert x_{0} \Vert + \int_{0}^{t}\Vert \dot{x}(s) \Vert\,\mathrm{d}s $ and noting that for a.e. $t\in [0,T]$ we get
$ 1+\Vert x(t) \Vert\leq \rho(t)$ give us by \eqref{e:rho} the estimate
\begin{equation*}
\begin{aligned}
\dot{\rho}(t)\leq m_{F}+\beta\int_{0}^{t}\rho(s)\,\mathrm{d}s\;\mbox{ for a.e. }\;t\in [0,T].
\end{aligned}
\end{equation*}
Applying the extended Gronwall's inequality from Lemma~\ref{granwal:cont} with $ \rho(\cdot)$ yields
\begin{equation*}
\begin{aligned}
\rho(t)&\le\big(1+\Vert x_{0} \Vert\big)\exp\Big(\int_{0}^{t}(\beta +1)\,\mathrm{d}\tau\Big)+\int_{0}^{t}m_{F}\exp\Big(\int_{s}^{t}(\beta +1)\,\mathrm{d}\tau\Big)\,\mathrm{d}s\\
&=(1+\Vert x_{0} \Vert)\exp\left(t(\beta+1)\right)+\dfrac{m_{F}}{\beta+1}\exp\left(t(\beta+1)\right).
\end{aligned}
\end{equation*}
This ensures the fulfillment of the estimate
\begin{equation*}
1+\Vert x(t)\Vert\le\Big(1+\Vert x_{0} \Vert +\dfrac{m_{F}}{\beta+1}\Big)\exp\big(T(\beta+1)\big),
\end{equation*}
and hence justifies the first inequality in \eqref{bound:solution}.
Consequently, we get
\begin{equation*}
\begin{aligned}
\Vert \dot{x}(t)\Vert\leq m_{F}+\beta T\,M_{1}\;\mbox{ for a.e. }\;t\in[0,T].
\end{aligned}
\end{equation*}
which verifies the second estimate in \eqref{bound:solution} and thus completes the proof of the theorem.$\h$\vspace*{-0.2in}

\section{Discrete Approximations of Integro-Differential Dynamics}\label{sec:discrete}
\setcounter{equation}{0}\vspace*{-0.1in}

In this section, we obtain one of the most significant results of the paper, which is important not only for deriving necessary optimality conditions for integro-differential systems, but certainly has its own value from the viewpoints of stability and numerical analysis of IDIs. Here we show that, under the imposed assumptions on the IDI data, any feasible solution to \eqref{inclusion} can be constructively approximated in the $W^{1,2}$-norm topology by piecewise-linear extensions on $[0,T]$ of feasible solutions to finite-dimensional inclusions with discrete time.

Considering for each fixed $k\in\N$ the discrete mesh \eqref{mesh} on $[0,T]$, we get the following strong $W^{1,2}$-approximation result for {\em any} feasible solution to \eqref{inclusion}.\vspace*{-0.1in}

\begin{theorem}[\bf strong discrete approximations of feasible solutions to IDIs]\label{strong approximation}
Let $ \bar{x}(\cdot) $ be a solution to the integro-differential inclusion $ (D_{F,g}) $ under assumptions $(\mathcal{H}^{F})$ and $(\mathcal{H}^{g})$, and let $\Delta_k$ be the discrete mesh on $[0,T]$ defined in \eqref{mesh}. Then for any fixed $k\in\N$, there exist a piecewise-linear function $x^k(\cdot)$ and a piecewise-constant function $y^k(\cdot)$ on $ [0,T] $ such that
$(x^{k}(0),y^{k}(0))=(x_{0},0)$ and
\begin{equation}\label{disc1}
x^{k}(t)=x^{k}(t^{k}_{j})+(t-t^{k}_{j})\upsilon^{k}_{j},\quad y^{k}(t)=w^{k}_{j}\;\mbox{ for }\;t\in(t^{k}_{j},t^{k}_{j+1}],\;j=0,\ldots,k-1,
\end{equation}
where the discrete velocities $\upsilon^{k}_{j}$ and $w^{k}_{j}$ satisfy the conditions
\begin{equation}\label{disc2}
\upsilon^{k}_{j}\in F\big(t^{k}_{j},x^{k}(t^{k}_{j})\big)+w^{k}_{j},\;\;j=0,\ldots,k-1,
\end{equation}
\begin{equation}\label{disc3}
w^{k}_{j}:=\dfrac{1}{h^{k}_{j}}\int\limits_{t^{k}_{j}}^{t^{k}_{j+1}}\Big\{\sum\limits_{i=0}^{j-1}\int\limits_{t^{k}_{i}}^{t^{k}_{i+1}}g\big(t,s,x^{k}(t^{k}_{i})\big)\,
\mathrm{d}s+\int\limits_{t^{k}_{j}}^{t}g\big(t,s,x^{k}(t^{k}_{j})\big)\,\mathrm{d}s\Big\}\,\mathrm{d}t,\;\;j=0,\ldots,k-1.
\end{equation}
Moreover, we have the  convergence
$x^{k}(\cdot)\to \bar{x}(\cdot)$ in the $W^{1,2}$-norm topology and $ y^{k}(\cdot)\to\oy(\cdot)$ in the $L^{2}$-norm topology on $ [0,T]$ as $k\to\infty$, where $\oy(\cdot)$ is defined by
\begin{equation}\label{e2:bkj}
\bar{y}(t):=\int_{0}^{t}g\big(t,s,\bar{x}(s)\big)\,\mathrm{d}s\;\mbox{ for all }\;t\in[0,T].
\end{equation}
\end{theorem}\vspace*{-0.1in}
{\bf Proof}. Recalling that step functions are dense in $ L^{2}([0,T];\mathbb{R}^{n})$ gives us a sequence $\{a^{k}(\cdot)\} $  with
\begin{equation*}
\xi_{k}^{2}:=T\int_{0}^{T}\Vert a^{k}(t)-\dot{\bar{x}}(t)\Vert^{2}\,\mathrm{d}t\rightarrow 0 \quad \text{as} \quad k\rightarrow \infty.
\end{equation*}
For each $ k\in \mathbb{N} $, there exists a partition $ \Delta_{k} $ of the interval $ [0,T] $ in \eqref{mesh} for which
the step functions $ \{a^{k}(\cdot)\} $ are constant on $ [t^{k}_{j},t^{k}_{j+1}) $ as $ j=0,\ldots, k-1 $.
For simplicity, we use the notation $t_j:=t^k_j$ for the mesh points as $j=0,\ldots,k$ and $k\in\N$.
\vspace{0.3cm}\\
Let us define the piecewise-linear functions $ u^{k}(\cdot)$ on $[0,T]$ by
\begin{equation}\label{uk}
u^{k}(0):=x_{0},\;\;u^{k}(t):=u^{k}(t_{j})+(t-t_{j})a^{k}(t_{j})\;\mbox{ for }\;t\in[t^{k}_{j},t^{k}_{j+1}],\;j=0,\ldots,k-1,
\end{equation}
and observe that these functions are piecewise-linear extensions of $ u^{k}(t_{j})$ to the continuous-time interval $[0,T]$, admit the representations
\begin{equation*}
u^{k}(t)=x_{0}+\int_{0}^{t}a^{k}(s)\,\mathrm{d}s\;\mbox{ whenever }\;t\in[0,T],
\end{equation*}
and for all $k\in\N$ satisfy the estimates
\begin{equation}\label{e:ineq}
\Vert u^{k}(t)-\bar{x}(t)\Vert\leq\int_{0}^{T}\Vert a^{k}(t)-\dot{\bar{x}}(t)\Vert\,\mathrm{d}t\leq\xi_{k},\quad t\in[0,T].
\end{equation}
Next we construct the piecewise-constant mappings $ b^{k}: [0,T]\rightarrow\mathbb{R}^{n} $ by
\begin{equation}\label{e:bk}
b^{k}(0):=0\quad\text{and}\quad b^{k}(t):=b^{k}_{j},\;\mbox{ if }\;t\in(t^{k}_{j},t^{k}_{j+1}],\;\;j=0,\ldots,k-1,
\end{equation}
through their values at the mesh points
\begin{equation}\label{e:bkj}
b^{k}_{j}:=\dfrac{1}{h^{k}_{j}}\int\limits_{t^{k}_{j}}^{t^{k}_{j+1}}\Big\{\sum\limits_{i=0}^{j-1}\int\limits_{t^{k}_{i}}^{t^{k}_{i+1}}g\big(t,s,u^{k}(t^{k}_{i})\big)\,\mathrm{d}s
+\int\limits_{t^{k}_{j}}^{t}g\big(t,s,u^{k}(t^{k}_{j})\big)\,\mathrm{d}s\Big\}\,\mathrm{d}t,\;\;j=0,\ldots,k-1,
\end{equation}
and then define the piecewise-constant functions $\theta^{k}(\cdot):[0,T]\longrightarrow[0,T]$ by
\begin{equation}\label{e:uk}
\theta^{k}(0):=0\quad\text{and}\quad
\theta^{k}(t):=t_{j}\;\text{ if }\;t\in (t_{j},t_{j+1}].
\end{equation}
Using the constructions in \eqref{e:bk} and \eqref{e:bkj}, we get for any $ j=0,\ldots,k-1 $ and $k\in\N$ that
\begin{equation}\label{e2:bk}
b^{k}(t)=\dfrac{1}{h^{k}_{j}}\int\limits_{t^{k}_{j}}^{t^{k}_{j+1}}\Big\{\int_{0}^{t}g\big(t,s,u^{k}(\theta^{k}(s))\big)\,\mathrm{d}s\Big\}\,\mathrm{d}t.
\end{equation}

Our next goal is to prove the $ L^{2}$-strong convergence of the sequence $ \{b^{k}(\cdot)\} $ to $ \bar{y}(\cdot) $,  where $\oy(\cdot)$ is defined in \eqref{e2:bkj}. Since $ \bar{x}(\cdot)\in W^{1,2}([0,T],\mathbb{R}^{n})$, we can choose a constant $M>0$ such that $\lVert \bar{x}(t)\rVert \leq M $. It follows from the integrable linear growth conditions $ (\mathcal{H}^{g}_{1}) $ that
\begin{equation*}
\Vert\bar{y}(t)\Vert\leq \beta\int_{0}^{t}(1+\Vert\bar{x}(s)\Vert)\,\mathrm{d}s\leq T\beta(1+M).
\end{equation*}
Employing the above constructions and the assumptions in $ (\mathcal{H}^{g}_{2}) $ tells us that
\begin{equation}\label{e:y}
\begin{aligned}
\left\Vert\bar{y}(t)-\int_{0}^{t}g\left(t,s,u^{k}(\theta^{k}(s))\right)\,\mathrm{d}s\right\Vert&\leq\int_{0}^{t}\Vert g\left(t,s,\bar{x}(s)\right)-g\left(t,s,u^{k}(\theta^{k}(s))\right)\Vert\,\mathrm{d}s\\
&\leq\alpha\int_{0}^{T}\Vert \bar{x}(t)-u^{k}(\theta^{k}(t))\Vert\,\mathrm{d}t\\
&\leq\alpha\int_{0}^{T}\Vert \bar{x}(t)-u^{k}(t)\Vert\,\mathrm{d}t+\alpha\int_{0}^{T}\Vert u^{k}(t)-u^{k}(\theta^{k}(t))\Vert\,\mathrm{d}t
\end{aligned}
\end{equation}
for all $ t\in[0,T]$. Moreover, it follows from \eqref{e:uk} that
\begin{equation*}
\begin{aligned}
\Big\Vert u^{k}(t)-u^{k}\big(\theta^{k}(t)\big)\Big\Vert &=\Big\Vert\int_{\theta^{k}(t)}^{t}a^{k}(\tau)\,\mathrm{d}\tau\Big\Vert\leq \int_{\theta^{k}(t)}^{t}\big\Vert a^{k}(\tau)\big\Vert
\,\mathrm{d}\tau\\
&\leq \int_{0}^{T}\Vert a^{k}(t)-\dot{\bar{x}}(t)\Vert\,\mathrm{d}t+\int_{\theta ^{k}(t)}^{t}\big\Vert \dot{\bar{x}}(\tau)\big\Vert
\,\mathrm{d}\tau.
\end{aligned}
\end{equation*}
Substituting the latter estimates together with \eqref{e:ineq} into \eqref{e:y} yields
\begin{equation*}
\left\Vert\bar{y}(t)-\int_{0}^{t}g\left(t,s,u^{k}(\theta^{k}(s))\right)\,\mathrm{d}s\right\Vert\leq \alpha(T+1)\xi_{k} +\alpha\int_{0}^{T}\left\{\int_{\theta ^{k}(t)}^{t}\left\Vert \dot{\bar{x}}(\tau)\right\Vert
\,\mathrm{d}\tau\right\}\,\mathrm{d}t,
\end{equation*}
which being combined with \eqref{e2:bk} leads us to the inequalities
\begin{equation*}
\begin{aligned}
&\Vert b^{k}(t)-\bar{y}(t)\Vert\leq\Big\Vert b^{k}(t)-\dfrac{1}{h^{k}_{j}}\int\limits_{t^{k}_{j}}^{t^{k}_{j+1}}\bar{y}(t)\,\mathrm{d}t\Big\Vert + \Big\Vert\dfrac{1}{h^{k}_{j}}\int\limits_{t^{k}_{j}}^{t^{k}_{j+1}}\bar{y}(t)\,\mathrm{d}t-\bar{y}(t)\Big\Vert\\
&\leq\alpha(T+1)\xi_{k}+\alpha\int_{0}^{T}\Big\{\int_{\theta ^{k}(t)}^{t}\Big\Vert \dot{\bar{x}}(\tau)\Big\Vert
\,\mathrm{d}\tau\Big\}\,\mathrm{d}t +\Big\Vert\dfrac{1}{h^{k}_{j}}\int\limits_{t^{k}_{j}}^{t^{k}_{j+1}}\bar{y}(t)\,\mathrm{d}t-\bar{y}(t)\Big\Vert.
\end{aligned}
\end{equation*}
Since $\dot{\bar{x}}(\cdot)\in L^{1}([0,T],\mathbb{R}^{n})$ and for each $t\in [0,T]$ we have $\theta^{k}(t)\longrightarrow t$, it follows that $ \lim\limits_{k\rightarrow +\infty}\int_{\theta ^{k}(t)}^{t}\left\Vert \dot{\bar{x}}(\tau)\right\Vert\,\mathrm{d}\tau=0 $. Remembering that $ \xi_{k}\rightarrow 0 $ and applying Lemma~\ref{mu:lemma} justify the
strong convergence of $\{b^{k}(\cdot)\}$ to $ \bar{y}(\cdot) $ in the norm topology of $L^{2}([0,T],\mathbb{R}^{n})$ as $k\to\infty$.

Now we use the discrete functions $ a^{k}(t_{j})$ and $ b^{k}(t_{j})$ to construct $ x^{k}(t_{j})$, $ v^{k}(t_{j})$, and $w^{k}(t_{j})$ at the mesh points of $\Delta_k$ such that their corresponding extensions on $[0,T]$ satisfy the relationships in \eqref{disc1}--\eqref{disc3} and exhibit the claimed convergence properties. To furnish this, let us proceed recurrently by the following {\em projection algorithm}, which is clearly implementable due to the imposed assumptions:
\begin{equation}\label{prox}
\left\{
\begin{array}{l}
x^{k}(0)=x_{0},\;\;x^{k}(t_{j+1})=x^{k}(t_{j})+h^{k}_{j}\upsilon^{k}_{j},\;\;j=0,\ldots,k-1,\vspace*{0.3cm}\\
\upsilon^{k}_{j}\in F\big(t_{j},x^{k}(t_{j})\big)+w^{k}_{j},\;\;j=0,\ldots,k-1,\vspace*{0.3cm}\\
w^{k}_{j}=\disp\dfrac{1}{h^{k}_{j}}\int\limits_{t^{k}_{j}}^{t^{k}_{j+1}}\Big\{\sum\limits_{i=0}^{j-1}\int\limits_{t^{k}_{i}}^{t^{k}_{i+1}}g\big(t,s,x^{k}(t^{k}_{i})\big)\,\mathrm{d}s
+\int\limits_{t^{k}_{j}}^{t}g\big(t,s,x^{k}(t^{k}_{j})\big)\,\mathrm{d}s\Big\}\,\mathrm{d}t,\;\;j=0,\ldots,k-1,\vspace*{0.3cm}\\
\Vert(\upsilon^{k}_{j}-a^{k}(t_{j}))+(b^{k}(t_{j})-w^{k}_{j}) \Vert=\mathrm{dist}\big(a^{k}(t_{j})-b^{k}(t_{j});F(t_{j},x^{k}(t_{j})\big),\;\;j=0,\ldots,k-1.
\end{array}\right.
\end{equation}
It is easy to see the constructions in \eqref{prox} with $y^k(\cdot)$ defined as in \eqref{disc1} satisfy all the conditions in \eqref{disc1}--\eqref{disc3}. We need to verify the convergence properties claimed in the theorem.

To this end, fix any $ 0\leq j\leq k-1 $ and $ s\in[t_{j},t_{j+1}] $ where $\dot{\bar{x}}(s)$ exists. Then we have
\begin{equation}\label{e:vk}
\begin{aligned}
\lVert\upsilon^{k}_{j}-\dot{\bar{x}}(s) \rVert&\leq \Vert \upsilon^{k}_{j}-a^{k}(s)+b^{k}(s)-w^{k}_{j} \Vert + \Vert a^{k}(s)-\dot{\bar{x}}(s) \Vert + \Vert w^{k}_{j}-b^{k}(s)\Vert.
\end{aligned}
\end{equation}
On the other hand, it follows from \eqref{tau} and the assumptions in $ (\mathcal{H}^{F}_{2}) $ that
\begin{equation*}
\begin{aligned}
&\Vert \upsilon^{k}_{j}-a^{k}(s)+b^{k}(s)-w^{k}_{j} \Vert =\mathrm{dist}\left(a^{k}(s)-b^{k}(s),F(t_{j},x^{k}(t_{j}))\right)\\
&\leq \mathrm{dist}\left(a^{k}(s)-b^{k}(s),F(t_{j},u^{k}(t_{j}))\right)+l_{F}\Vert x^{k}(t_{j})-u^{k}(t_{j}) \Vert\\
&\leq \mathrm{dist}\left(a^{k}(s)-b^{k}(s),F(s,u^{k}(t_{j}))\right)+l_{F}\Vert x^{k}(t_{j})-u^{k}(t_{j}) \Vert +\tau(F,h_{k})\\
&\leq \mathrm{dist}\left(a^{k}(s)-b^{k}(s),F(s,\bar{x}(s))\right)+l_{F}(\Vert u^{k}(t_{j})-u^{k}(s)\Vert + \Vert u^{k}(s)-\bar{x}(s) \Vert)\\
&+l_{F}\Vert x^{k}(t_{j})-u^{k}(t_{j})\Vert +\tau(F,h_{k}),
\end{aligned}
\end{equation*}
The definition of $ u^{k}(\cdot)$ in \eqref{uk} and the feasibility of $ \bar{x}(\cdot)$ in \eqref{inclusion} imply that
\begin{equation*}
\begin{aligned}
&\Vert \upsilon^{k}_{j}-a^{k}(s)+b^{k}(s)-w^{k}_{j} \Vert\\
&\leq \Vert a^{k}(s)-\dot{\bar{x}}(s) \Vert + \Vert b^{k}(s)-\bar{y}(s) \Vert +l_{F}\Vert u^{k}(s)-\bar{x}(s)\Vert + l_{F}\Vert x^{k}(t_{j})-u^{k}(t_{j}) \Vert\\
& +l_{F}(s-t_{j})\Vert a^{k}(s) \Vert +\tau(F,h_{k}).
\end{aligned}
\end{equation*}
Observe in addition due to the constructions above and the assumptions in $ (\mathcal{H}^{g}_{2}) $ that
\begin{equation*}
\begin{aligned}
\Vert w^{k}_{j}-b^{k}(s)\Vert&\leq \alpha\sum\limits_{i=0}^{j}h^{k}_{i}\Vert x^{k}(t_{i})-u^{k}(t_{i}) \Vert +\dfrac{\alpha}{h^{k}_{j}}\int\limits_{t_{j}}^{t{j+1}}(s-t_{j})\Vert x^{k}(t_{j})-u^{k}(t_{j})\Vert\,\mathrm{d}s\\
&\leq\alpha\sum\limits_{i=0}^{j}h^{k}_{i}\Vert x^{k}(t_{i})-u^{k}(t_{i}) \Vert +(\alpha h^{k}_{j}/2)\Vert x^{k}(t_{j})-u^{k}(t_{j})\Vert.
\end{aligned}
\end{equation*}
Substituting the latter estimates into \eqref{e:vk} gives us
\begin{equation*}
\begin{aligned}
\lVert\upsilon^{k}_{j}-\dot{\bar{x}}(s) \rVert &\leq 2\Vert a^{k}(s)-\dot{\bar{x}}(s) \Vert + \Vert b^{k}(s)-\bar{y}(s) \Vert +l_{F}\Vert u^{k}(s)-\bar{x}(s)\Vert+ (l_{F}+\alpha h^{k}_{j}/2)\Vert x^{k}(t_{j})-u^{k}(t_{j}) \Vert\\
& +\alpha\sum\limits_{i=0}^{j}h^{k}_{i}\Vert x^{k}(t_{i})-u^{k}(t_{i}) \Vert  + l_{F}(s-t_{j})\Vert a^{k}(s) \Vert +\tau(F,h_{k}).
\end{aligned}
\end{equation*}
Then we deduce from \eqref{e:ineq} the inequalities
\begin{equation*}
\begin{aligned}
\Vert x^{k}(t_{j})-u^{k}(t_{j}) \Vert\leq \Vert x^{k}(t_{j})-\bar{x}(t_{j}) \Vert + \Vert \bar{x}(t_{j})-u^{k}(t_{j}) \Vert\leq \Vert x^{k}(t_{j})-\bar{x}(t_{j}) \Vert + \xi_{k},
\end{aligned}
\end{equation*}
and therefore arrive at the estimates
\begin{equation}\label{e:wk-bk}
\begin{aligned}
\Vert w^{k}_{j}-b^{k}(s)\Vert &\leq\alpha\sum\limits_{i=0}^{j}h^{k}_{i}\Vert x^{k}(t_{i})-\bar{x}(t_{i}) \Vert +(\alpha h^{k}_{j}/2)\Vert x^{k}(t_{j})-\bar{x}(t_{j})\Vert +(\alpha\,T + \alpha h^{k}_{j}/2)\xi_{k},
\end{aligned}
\end{equation}
\begin{equation}\label{e:vk2}
\begin{aligned}
\lVert\upsilon^{k}_{j}-\dot{\bar{x}}(s)\rVert &\leq 2\Vert a^{k}(s)-\dot{\bar{x}}(s) \Vert + \Vert b^{k}(s)-\bar{y}(s) \Vert+l_{F}\Vert u^{k}(s)-\bar{x}(s)\Vert\\
&+ (l_{F}+\alpha h^{k}_{j}/2)\Vert x^{k}(t_{j})-\bar{x}(t_{j}) \Vert+\alpha\sum\limits_{i=0}^{j}h^{k}_{i}\Vert x^{k}(t_{i})-\bar{x}(t_{i}) \Vert\\
& + l_{F}(s-t_{j})\Vert a^{k}(s) \Vert +(l_{F}+\alpha\,T+\alpha h^{k}_{j}/2))\xi_{k}+\tau(F,h_{k}).
\end{aligned}
\end{equation}
Since the set of $s$ where $\dot\ox(s)$ exists is of full measure on $[0,T]$, we denote
\begin{equation*}
c_{k}(s):=2\Vert a^{k}(s)-\dot{\bar{x}}(s) \Vert + \Vert b^{k}(s)-\bar{y}(s) \Vert + l_{F}(s-t_{j})\Vert a^{k}(s) \Vert +(2l_{F}+\alpha\,T+\alpha h^{k}_{j}/2)\xi_{k}+\tau(F,h_{k}),
\end{equation*}
get  $\lim\limits_{k\to \infty}\int_{0}^{T}c^{2}_{k}(s)\,\mathrm{d}s=0$, and derive from the above relationships that
\begin{equation*}
\begin{aligned}
\lVert x^{k}(t_{j+1})&-\bar{x}(t_{j+1}) \rVert =\Big\lVert  x^{k}(t_{j})+h_{k}v^{k}_{j}-\bar{x}(t_{j})-\int\limits_{t_{j}}^{t_{j+1}}\dot{\bar{x}}(s)\,\,\mathrm{d}s \Big\rVert\leq\lVert x^{k}(t_{j})-\bar{x}(t_{j})
\rVert + \int\limits_{t_{j}}^{t_{j+1}}\lVert v^{k}_{j}-\dot{\bar{x}}(s)\rVert\,\mathrm{d}s\\
&\leq(1+h_{k}(l_{F}+3\alpha h_{k}/2)) \lVert x^{k}(t_{j})-\bar{x}(t_{j}) \rVert  +\alpha\,h_{k}^{2}\sum\limits_{i=0}^{j-1}\lVert x^{k}(t_{i})-\bar{x}(t_{i}) \rVert +\int\limits_{t_{j}}^{t_{j+1}}c_{k}(s)\,\mathrm{d}s.
\end{aligned}
\end{equation*}
Applying now the discrete Gronwall's lemma from Proposition~\ref{granwal} with the given parameters
\begin{equation*}
e_{j}=\lVert x^{k}(t_{j})-\bar{x}(t_{j}) \rVert,\quad\rho_{j}=\alpha\,h_{k}^{2},\quad\gamma_{j}=h_{k}(l_{F}+3\alpha h_{k}/2),\;\;\mbox{ and }\;\;\sigma_{j}=\int\limits_{t_{j}}^{t_{j+1}}c_{k}(s)\,\mathrm{d}s
\end{equation*}
yields for all $j=0,\ldots,k-1$ the estimates
\begin{equation}\label{e:xk2}
\begin{aligned}
\lVert x^{k}(t_{j})-\bar{x}(t_{j}) \rVert &\leq \int_{0}^{T}c_{k}(s)\,\mathrm{d}s\,\exp\Big(\alpha\dfrac{j(j-1)h_{k}^{2}}{2}+jh_{k}(l_{F}+3\alpha h_{k}/2)\Big)\\
& \leq \int_{0}^{T}c_{k}(s)\,\mathrm{d}s\,\exp\left(\alpha\,T^{2}/2+T(l_{F}+3\alpha h_{k}/2)\right):=\zeta_{k},
\end{aligned}
\end{equation}
where $ \lim\limits_{k\to\infty}\zeta_{k}=0 $. Hence we deduce from \eqref{e:vk2} that
\begin{equation}\label{e:vk3}
\lVert\upsilon^{k}_{j}-\dot{\bar{x}}(s) \rVert \leq c_{k}(s)+(l_{F}+2\alpha\,T+\alpha h_{k}/2)\zeta_{k}\;\mbox{ whenever }\;s\in[t_{j},t_{j+1}]
\end{equation}
on a full measure set. Employing the obtained conditions together with \eqref{e:xk2} tells us that
\begin{equation*}
\begin{aligned}
\lVert x^{k}(t)-\bar{x}(t) \rVert&=\Big\lVert x^{k}(t_{j})+(t-t_{j})\upsilon^{k}_{j}-\bar{x}(t_{j})-\int\limits_{t_{j}}^{t}\dot{\bar{x}}(s)\,\mathrm{d}s\Big\rVert\leq \lVert x^{k}(t_{j})-\bar{x}(t_{j}) \rVert + \int\limits_{t_{j}}^{t_{j+1}}\lVert v^{k}_{j}-\dot{\bar{x}}(s)\rVert\,\mathrm{d}s\\
&\leq \zeta_{k}+\int_{0}^{T}c_{k}(s)\,\mathrm{d}s + h_k(l_{F}+2\alpha\,T+\alpha h_{k}/2)\zeta_{k}\;\mbox{ for all }\;t\in[0,T],
\end{aligned}
\end{equation*}
which readily verifies the uniform convergence of the sequence $\{x^{k}(\cdot)\}$ to $\bar{x}(\cdot) $ as $k\to \infty $.

To proceed further, we deduce from \eqref{e:vk3} that
\begin{equation}\label{beta}
\begin{aligned}
\int_{0}^{T}\lVert \dot{x}^{k}(t)-\dot{\bar{x}}(t) \rVert^{2}\,\mathrm{d}t
&=\sum\limits_{j=0}^{k-1}\int\limits_{t_{j}}^{t_{j+1}}\lVert \upsilon^{k}_{j}-\dot{\bar{x}}(t) \rVert^{2}\,\mathrm{d}t\\
&\le\int_{0}^{T}c^{2}_{k}(t)\,\mathrm{d}t +T(l_{F}+2\alpha\,T+\alpha h_{k}/2)^{2}\zeta^{2}_{k}:=\beta_{k}.
\end{aligned}
\end{equation}
Remembering that $\zeta_{k}\to 0$ and $ \int_{0}^{T}c^{2}_{k}(t)\,\mathrm{d}t\to 0 $, this justifies the strong convergence of the sequence $\{\dot{x}^{k}(\cdot)\}$ to $ \dot{\bar{x}}(\cdot) $ in the norm topology of $L^{2}([0,T],\mathbb{R}^{n})$ and therefore the claimed $W^{1,2}$-norm convergence of the extended discrete trajectories $\{x^{k}(\cdot)\}$ to the feasible one $\bar{x}(\cdot)$  for IDE \eqref{inclusion} as $k\to \infty $.

To complete the proof of the theorem, it remains to verify the $ L^{2}$-strong convergence of the sequence $\{y^{k}(\cdot)\}$ to $\bar{y}(\cdot)$ as $k\to\infty$, where $ y^{k}(t)=w^{k}_{j} $ for all $ t\in(t^{k}_{j},t^{k}_{j+1}] $ and $ j=0,\ldots,k-1 $. To this end, observe first by \eqref{e:wk-bk} and \eqref{e:xk2} that for any $ s\in[t_{j},t_{j+1}]$ we have
\begin{equation*}
\begin{aligned}
\Vert y^{k}(s)-b^{k}(s) \Vert &\leq (\alpha\,T+\alpha\,h_{k}/2)\zeta_{k}+(\alpha\,T + \alpha h^{k}_{j}/2)\xi_{k}\\
&\leq(\alpha\,T + \alpha h_{k}/2)(\zeta_{k}+\xi_{k}).
\end{aligned}
\end{equation*}
The latter ensures the fulfillment of the relationships
\begin{equation*}
\begin{aligned}
\int_{0}^{T}\Vert y^{k}(s)-\bar{y}(s)\Vert^{2}\,\mathrm{d}s&=\sum\limits_{j=0}^{k-1}\int\limits_{t_{j}}^{t_{j+1}}\Vert y^{k}(s)-\bar{y}(s)\Vert^{2}\mathrm{d}s\leq \sum\limits_{j=0}^{k-1}\int\limits_{t_{j}}^{t_{j+1}}\Vert y^{k}(s)-b^{k}(s)\Vert^{2}\,\mathrm{d}s + \int_{0}^{T}\Vert b^{k}(s)-\bar{y}(s)\Vert^{2}\,\mathrm{d}s \\
&\leq T(\alpha\,T + \alpha h_{k}/2)^{2}(\zeta_{k}+\xi_{k})^{2}+\int_{0}^{T}\Vert b^{k}(s)-\bar{y}(s)\Vert^{2}\,\mathrm{d}s,
\end{aligned}
\end{equation*}
which justify therefore the desired strong convergence of $ y^{k}(\cdot) $ to $ \bar{y}(\cdot) $ in $ L^{2}([0,T],\mathbb{R}^{n})$. $\h$\vspace*{-0.2in}

\section{Intermediate Local Minimizers and Relaxation Stability}\label{sec:loc-min}
\setcounter{equation}{0}\vspace*{-0.1in}

This section concerns some fundamental properties of dynamic optimization problems for integro-differential systems over continuous-time intervals that are important for the method of discrete approximations and allows us, in particular, to select an appropriate class of local minimizers in $(P)$ for which we derive necessary optimality conditions.

It has been well recognized (starting with Weierstrass at the end of 19th century) that continuous-time problems of the calculus of variations, optimal control, and dynamic optimization with nonconvex velocities may not admit optimal solutions. On the other hand, Bogolyubov and Young proved at the beginning of 1930s that the {\em convexification} of classical variational problems with respect to velocity variables produced an extended/relaxed problem admitting an optimal solution, which can be uniformly approximated by a minimizing sequence of feasible solutions to the original problem. This implies, in particular, that optimal values of the cost functionals in the original and relaxed problems agree, meaning that the property of {\em  relaxation stability} holds. This line of research has been further developed for various classes of continuous-time dynamic optimization problems in many publications; see, e.g., \cite{Blasi,Boris} for more details, references, and historical remarks.

Following this direction for the integro-differential optimization problem $(P)$, we define the convexified multifunction
$\mathrm{co}\,F:[0,T]\times\mathbb{R}^{n} \rightrightarrows \mathbb{R}^{n}$ by
\begin{equation*}
(\mathrm{co}\,F)(t,x):=\mathrm{co}\,F(t,x) \quad\text{for all}\;\; t\in [0,T]\;\mbox{ and }\;x \in \mathbb{R}^{n},
\end{equation*}
and then consider, along with $(D_{\mathrm{co}\,F,g})$ in \eqref{inclusion}, the {\em relaxed IDI} formulated as
\begin{equation}\label{co:inclusion}
(D_{\mathrm{co}\,F,g}) :\;\;\left\{
\begin{array}{l}
\dot{x}(t)\in \mathrm{co}\,F\big(t,x(t)\big)+\displaystyle\int_{0}^{t}g\big(t,s,x(s)\big)\,\mathrm{d}s\quad \text{a.e.}\;\;t\in [0,T],\vspace*{0.3cm}\\
x(0)=x_0.
\end{array}\right.
\end{equation}
Accordingly, an arc $ x(\cdot)\in W^{1,2}$ satisfying \eqref{co:inclusion} is called a {\em relaxed trajectory} for $(D_{F,g})$. It follows from the results of \cite{Gaouir} that under the imposed assumptions on $F$ and $g$, we can uniformly approximate any relaxed trajectory of IDI but those for \eqref{inclusion}. Note that neither the cost functional \eqref{Bolza-functional}, not the endpoint constraint \eqref{endpoint constraints} are involved in the proposition below.\vspace*{-0.1in}

\begin{proposition}[\bf approximation property for relaxed IDI trajectories]\label{app-traj} Any absolutely continuous trajectory of the convexified IDI \eqref{co:inclusion} can be represented as the limit in the space $\mathcal{C}([0,T],\mathbb{R}^n)$ of a sequence $ x_{k}(\cdot) $ of absolutely continuous solutions to the original IDI in \eqref{inclusion}.
\end{proposition}\vspace*{-0.1in}

Now we define the set of {\em relaxed feasible solutions} satisfying the convexified IDI \eqref{co:inclusion} subject to the endpoint constraint \eqref{endpoint constraints} as follows:
\begin{equation}\label{feas-rela}
\mathcal{R}_{\mathrm{co}\,F,g}(x_0,\Omega):=\big\{x\in\mathcal{C}([0,T],\mathbb{R}^{n})\;\big|\;x(\cdot)\;\;\text{is a solution of}\;(D_{\mathrm{co}\,F,g})\;\;\text{with}\;\;x(T)\in\Omega\big\}.
\end{equation}
The next theorem establishes the compactness of the set $ \mathcal{R}_{\mathrm{co}\,F,g}(x_0,\Omega)$ in the space $\mathcal{C}([0,T],\mathbb{R}^n)$.\vspace*{-0.1in}

\begin{theorem}[\bf compactness of feasible solutions to relaxed constrained IDIs]\label{comp}
Suppose that the assumptions in $(\mathcal{H}^{F})$, $(\mathcal{H}^{g})$, and $ (\mathcal{H}^{\Omega})$ are satisfied. Then the set $ \mathcal{R}_{\mathrm{co}\,F,g}(x_0,\Omega)$ in \eqref{feas-rela} is compact with respect to the supremum norm topology.
\end{theorem}\vspace*{-0.1in}
{\bf Proof}. It is easy to observe that the assumptions in $(\mathcal{H}^{F})$ yield the fulfillment of their counterparts $(\mathcal{H}^{{\rm co}\,F})$ for the convexified multifunction $\mathrm{co}\,F$. Therefore, the existence results of Theorem~\ref{exist:solu} keep holding for $\mathrm{co}\,F$. This tells us that the set of relaxed feasible solutions $ \mathcal{R}_{\mathrm{co}\,F,g}(x_0,\Omega) $ is nonempty and bounded with respect to the supremum norm. Taking any sequence of feasible solutions $\{x^{k}(\cdot)\} $ to problem $(P) $, we aim at showing that there exists a solution $ x(\cdot) $ of $(D_{\mathrm{co}\,F,g})$ satisfying the constraints $ x(T)\in\Omega $ such that $ x^{k}(\cdot)\rightarrow x(\cdot)$ uniformly on $[0,T]$ along some subsequence. Indeed, it follows from \eqref{bound:solution} that the sequence $\{\dot{x}^{k}(\cdot)\}$ is bounded in $L^2([0,T],\mathbb{R}^{n})$. Remembering that bounded sets in reflexive spaces are weakly compact gives us some $\upsilon(\cdot)\in L^{2}([0,T],\mathbb{R}^{n})$ such that $\dot{x}_{k}(\cdot)\to\upsilon(\cdot)\in L^{2}([0,T],\mathbb{R}^{n})$ in the weak topology of this space along a subsequence (without relabeling). Define now the function
\begin{equation*}
x(t):=x_{0}+\int_{0}^{t}\upsilon(s)\,\mathrm{d}s\;\;\mbox{ for all }\;\;t\in[0,T],
\end{equation*}
and observe that $ x(0)=x_{0} $, $\dot{x}(t)=\upsilon(t)$ for a.e.\ $t\in[0,T]$, and thus $x(\cdot)\in W^{1,2}([0,T],\mathbb{R}^{n})$. We clearly have that $ \dot{x}^{k}(\cdot)\rightarrow \dot{x}(\cdot) $ weakly in $ L^{2}([0,T],\mathbb{R}^{n})$ and deduce from the weak convergence of derivatives that $ x^{k}(t) $ converges to $ x(t) $ for every $ t\in[0,T] $ as $ k\rightarrow\infty $.
Defining further
\begin{equation*}
y^{k}(t):=\int_{0}^{t}g\big(t,s,x^{k}(s)\big)\,\mathrm{d}s\;\;\text{and}\;\;y(t):=\int_{0}^{t}g\big(t,s,x(s)\big)\,\mathrm{d}s
\end{equation*}
allows us to find by \eqref{bound:solution} a constant $c>0$ such that for each $k\in\N$ we get $\lVert x^{k}(t)\rVert \leq \eta $ and $\lVert x(t)\rVert\le\eta $ as $t\in [0,T]$.
By $(\mathcal{H}_{2}^{g})$, there exists $\alpha>0 $ such that $ g(t,s,\cdot) $ is $ \alpha $-Lipschitz on $ \eta\B $. Then
\begin{equation}\label{cff}
\int_{0}^{T}\left\Vert y^{k}(t)-y(t)\right\Vert^{2}\,\mathrm{d}t \leq \int_{0}^{T} \alpha^{2}\left(\int_{0}^{t}\left\Vert x^{k}(s)-x(s) \right\Vert\,\mathrm{d}s\right)^{2}\,\mathrm{d}t.
\end{equation}
Note that for every $(t,s)\in Q_{\Delta} $ we have
\begin{equation*}
\alpha\int_{0}^{t}\left\Vert x^{k}(s)-x(s) \right\Vert\,\mathrm{d}s\leq 2\eta T\alpha.
\end{equation*}
If follows from \eqref{cff} and the Lebesgue dominated convergence theorem that
\[
y^{k}(\cdot)\rightarrow y(\cdot)\;\; \text{ strongly in }
\;\;L^{2}([0,T],\mathbb{R}^{n}),
\]
which implies in turn the convergence
\[
\zeta^{k}(\cdot):=\dot{x}^{k}(\cdot)-y^{k}(\cdot)\longrightarrow \zeta(\cdot):= \dot{x}(\cdot)-y(\cdot),
\]
weakly in $ L^{2}([0,T],\mathbb{R}^{n}) $. By Mazur's weak closure theorem, we find a sequence of convex combinations $ \sum\limits_{m=k}^{r(k)}\beta_{m,k}\zeta^{m}(\cdot) $ with $ \sum\limits_{m=k}^{r(k)}\beta_{m,k}=1 $ and $\beta_{m,k}\in[0,1] $ for all $ m,k $ that converges strongly in $ L^{2}([0,T],\mathbb{R}^{n}) $ to $ \zeta(\cdot)$. Extracting a subsequence tells us without loss of generality that $ \sum\limits_{m=k}^{r(k)}\beta_{m,k}\zeta^{m}(\cdot)$ converges almost everywhere on $ [0,T] $ to the mapping $ \zeta(\cdot) $. Remembering that $ x^{k}(\cdot) $ is sequence of feasible solutions to problem $ (P) $, we arrive at the inclusion
\begin{equation*}
\sum\limits_{m=k}^{r(k)}\beta_{m,k}\zeta^{m}(t)\in F\big(t,x^{k}(t)\big)\quad\text{ a.e.}\;\;t\in [0,T].
\end{equation*}
Observe furthermore by the assumptions in $ (\mathcal{H}^{F}_{2}) $ that
\begin{equation*}
\begin{aligned}
\mathrm{dist}\,\big(\zeta(t);\mathrm{co}\,F(t,x(t))\big)&=\mathrm{dist}\,\big(\zeta(t);\mathrm{co}\,F(t,x(t))\big)-\mathrm{dist}\,\Big(\sum\limits_{m=k}^{r(k)}\beta_{m,k}\zeta^{m}(t);
\mathrm{co}\,F\big(t,x^{k}(t)\big)\Big)\\
&\leq \Vert \zeta(t)-\sum\limits_{m=k}^{r(k)}\beta_{m,k}\zeta^{m}(t) \Vert + l_{F}\Vert x(t)-x^{k}(t) \Vert.
\end{aligned}
\end{equation*}
Passing to the limit as $k\to\infty$ in the above inequality and using the convergence $ \sum\limits_{m=k}^{r(k)}\beta_{m,k}\zeta^{m}(t)\rightarrow \zeta(t)$ for a.e. $ t\in [0,T] $ and $ x^{k}(t)\rightarrow x(t) $ for every $ t\in [0,T] $ ensure that
\begin{equation*}
\dot{x}(t)\in \mathrm{co}\,F\big(t,x(t)\big)+\int_{0}^{t}g\big(t,s,x(s)\big)\,\mathrm{d}s\quad\text{a.e.,}\;\;t\in [0,T].
\end{equation*}
Finally, we see that $ x(T)\in\Omega $ by the convergence $x^k(T)\to x(T)$ and the closedness of $ \Omega $ and thus confirm as claimed that the set $ \mathcal{R}_{\mathrm{co}\,F,g}(x_0,\Omega) $ is compact with respect to the supremum norm. $ \h $\vspace*{0.05in}

Note that our previous considerations in this section did not address the cost functional \eqref{Bolza-functional}. To include it into the relaxation process, we proceed as follows. Letting
\begin{equation*}
G(t,x):=F(t,x)+\int_{0}^{t}g\big(t,s,x(s)\big)\,\mathrm{d}s,
\end{equation*}
define the  {\em extended running cost} for $(P)$ by
\begin{equation*}
l_{G}(t,x,\upsilon):=l(t,x,\upsilon)+\delta_G(\upsilon)\;\mbox{ with }\;G:=G(t,x),
\end{equation*}
where the indicator function $\dd_G$ of the set $G$ is given by $\dd_G(\upsilon):=0$ if $\upsilon\in G$ and $\dd_G(\upsilon):=\infty$ otherwise. Then denote by $\hat{l}_{G}(t,x,\upsilon)$ the {\em biconjugate function} to $l_{G}(t,x,\upsilon)$ with respect to the velocity variable
\begin{equation*}
\hat{l}_{G}(t,x,\upsilon):=(l_{G})^{**}(t,x,\upsilon),
\end{equation*}
which is the largest proper, convex, and l.s.c.\ function with respect to $\upsilon$ majorized by $l_{G}$. The {\em relaxed problem} $(R)$ associated with the original problem $(P)$ for integro-differential inclusions is
\begin{equation}\label{R}
\text{minimize}\;\;\hat{J}[x]:=\varphi\big(x(T)\big)+\int_{0}^{T}\hat{l}_{G}\big(t,x(t),\dot{x}(t)\big)\,\mathrm{d}t,
\end{equation}
with $x(0)=x_0$ and $ x(T)\in \Omega $. Observe that if $ \hat{J}[x]<\infty $, then $ x(\cdot) $ satisfies the convexified integro-differential inclusion in \eqref{co:inclusion}. Arguing similarly to the usual relaxation procedure in problems of dynamic optimization (see, e.g., \cite[Theorem~1.1]{Blasi} and \cite[Theorem~6.11]{Boris}) with taking into account the above results of Proposition~\ref{app-traj} and Theorem~\ref{comp}, we arrive at the following approximation property for Bolza problems governed by integro-differential inclusions.\vspace*{-0.1in}

\begin{theorem}[\bf approximation property for IDI Bolza problems]\label{app-bolza}
Let $ x(\cdot) $ be a relaxed trajectory for the integro-differential inclusion $ (D_{F,g}) $, where $ F $ and $ g $ satisfies the assumptions in $ (\mathcal{H}^{F}) $ and $ (\mathcal{H}^{g})$, respectively. Suppose also that the integrand $ l $ in \eqref{Bolza-functional} is continuous in $ (x,\upsilon) $, measurable in $ t $, and uniformly bounded by a summable function near $ x(\cdot) $. Then there exists a sequence of the original trajectories $ x_{k}(\cdot) $ for $ (D_{F,g}) $ satisfying the conditions
\begin{equation*}
x_{k}(\cdot)\rightarrow x(\cdot)\quad\text{in}\quad\mathcal{C}([0,T],\mathbb{R}^{n}),
\end{equation*}
\begin{equation*}
\liminf\limits_{k\rightarrow\infty}\int_{0}^{T}l\big(t,x_{k}(t),\dot{x}_{k}(t)\big)\,\mathrm{d}t\leq\int_{0}^{T}\hat{l}_{G}\big(t,x(t),\dot{x}(t)\big)\,\mathrm{d}t.
\end{equation*}
\end{theorem}

Now we are ready to introduce the concepts of {\em local minimizers} in the integro-differential problem $(P)$ for which we derive necessary optimality conditions.\vspace*{-0.1in}

\begin{definition}[\bf intermediate local minimizers]\label{locmin} Consider the original problem $(P)$ in \eqref{Bolza-functional}--\eqref{endpoint constraints} and the associated relaxed problem given in \eqref{R}. Then we say that:\vspace*{-0.05in}

{\bf(i)} A feasible solution $\bar{x}(\cdot)$ to $(P)$ is an {\sc intermediate local minimizer} $($i.l.m.$)$ of this problem if $\bar{x}(\cdot)\in W^{1,2}([0,T];\mathbb{R}^{n})$ and there exists $\ve>0$ such that $J[\bar{x}]\le J[x]$ for any feasible solution $x(\cdot)$ to $(P)$ satisfying the $W^{1,2}$-localization conditions
\begin{equation}\label{loc}
\lVert \bar{x}(t)-x(t)\rVert <\varepsilon\;\;\text{for}\;\;t\in[0,T]\;\mbox{ and }\;\int_{0}^{T}\Vert \dot{\bar{x}}(t)-\dot{x}(t) \Vert^{2}\,\mathrm{d}t<\varepsilon.
\end{equation}\vspace*{-0.05in}

{\bf(ii)} A feasible solution $\bar{x}$ to $(P)$ is a {\sc relaxed intermediate local minimizer} $($r.i.l.m.$)$ of this problem if $J[\bar{x}]=\hat{J}[\bar{x}]$ and there exists $\varepsilon>0$ such that $J[\bar{x}]\le J[x]$ for any feasible solution $ x(\cdot) $ to $(P)$ satisfying  the localization conditions in \eqref{loc}.
\end{definition}\vspace*{-0.1in}

The i.l.m.\ notion from Definition~\ref{locmin}(i) is an integro-differential extension of the one first appeared in \cite{m95} for differential inclusions.  As discussed in \cite{m95,Boris}, this notion occupies an intermediate position between weak and strong minima for variational and control problems covering strong but not weak local minimizers. In the aforementioned references, the reader can find various examples showing that a local minimizer may be weak while not intermediate as well as intermediate but not strong; the latter may happen even for bounded, convex, and Lipschitzian differential inclusions.

Accordingly, the r.i.l.m.\ notion from Definition~\ref{locmin}(ii) is an integro-differential extension of the corresponding notion introduced in \cite{m95} for differential inclusions. We clearly have that any r.i.l.m.\ reduces to i.l.m.\ for problems with convex velocities. On the other hand, a remarkable feature of the continuous-time dynamical systems under consideration is their {\em hidden convexity}, which ensures that both notions of Definition~\ref{locmin} agree in general setting {\em without any convexity} assumptions. Indeed, Theorem~\ref{app-bolza} tells us that this is the case for problem $(P)$ in the absence of endpoint constraints. This also often happens for constrained nonconvex problems; see the discussions and references in \cite[Section~6.1.2]{Boris} for differential inclusions that can be readily extended to integro-differential ones.

As we see below, such a relaxation stability, being an internal property of continuous-time systems, allows us to implement the method of discrete approximations to deriving necessary optimality conditions for dynamic optimization problems governed by nonconvex constrained IDIs.\vspace*{-0.2in}

\section{Discrete Approximations of Local Minimizers}\label{sec:disc-opt}
\setcounter{equation}{0}\vspace*{-0.1in}

In this section, we construct a sequence of finite-dimensional dynamic optimization problems $(P_k)$ with discrete time such that piecewise-linear extensions of optimal solutions to $(P_k)$ on the continuous-time interval $[0,T]$ strongly converge in the
$W^{1,2}$-norm topology to a prescribed r.i.l.m.\ of $(P)$.

Given an r.i.l.m. $ \bar{x}(\cdot)$ for $(P)$ with $\ve>0$ from Definition~\ref{locmin}(ii), consider the mesh $\Delta_k$ defined in \eqref{mesh} and construct the sequence of discrete approximation problems  $(P_k)$, $k\in\N$, as follows:
\begin{equation}\label{disc-bolza}
\begin{aligned}
\text{minimize}\,\,\, J_{k}(x^{k})&:=\varphi(x^{k}_k)+h_{k}\sum\limits_{j=0}^{k-1}l\Big(t^{k}_{j},x^{k}_{j},\dfrac{x^{k}_{j+1}-x^{k}_{j}}{h_{k}}\Big)+\dfrac{1}{2}\sum\limits_{j=0}^{k-1}\int\limits_{t^{k}_{j}}^{t^{k}_{j+1}}\Big\lVert \dfrac{x^{k}_{j+1}-x^{k}_{j}}{h_{k}}-\dot{\bar{x}}(t) \Big\rVert^{2}\,\mathrm{d}t,
\end{aligned}
\end{equation}
over the collections $ (x_{0}^{k},\ldots,x^{k}_k)$ subject to the constraints
\begin{equation}\label{e:5.0}
\dfrac{x^{k}_{j+1}-x^{k}_{j}}{h^{k}_{j}}\in F(t^{k}_{j},x^k_j)+w^{k}_{j},\quad j=0,\ldots,k-1,
\end{equation}
with $ x^{k}_{0}=x_{0}$, where the vectors $ w^{k}_{j} $ are defined in \eqref{disc3}, and where
\begin{equation}\label{e:5.1}
x^{k}_{k}\in \Omega +\zeta_{k}\B\quad \text{with}\;\; \zeta_{k}\;\; \text{taken from}\;\; \eqref{e:xk2},
\end{equation}
\begin{equation}\label{e:5.2}
\lVert x^{k}_{j}-\bar{x}(t^{k}_{j}) \rVert\leq\dfrac{\varepsilon}{2},\,\,\,\,\, j=0,\ldots,k-1,\quad\mbox{and}
\end{equation}
\begin{equation}\label{e:5.3}
\begin{aligned}
&\sum\limits_{j=0}^{k-1}\int\limits_{t^{k}_{j}}^{t^{k}_{j+1}}\Big\lVert \dfrac{x^{k}_{j+1}-x^{k}_{j}}{h_{k}}-\dot{\bar{x}}(t)\Big\rVert^{2}\,\mathrm{d}t\leq\dfrac{\varepsilon}{2}.
\end{aligned}
\end{equation}

To proceed further, we need to make sure that problems $(P_k)$ admit optimal solutions.\vspace*{-0.1in}

\begin{proposition}[\bf existence of optimal solutions to discrete approximations]\label{discri:opti} Let $\bar{x}(\cdot)$ be an r.i.l.m.\ for problem $ (P)$ in \eqref{Bolza-functional}--\eqref{endpoint constraints}, and let the assumptions in $(\mathcal{H}^{F})$,
$(\mathcal{H}^{g})$, $(\mathcal{H}^{\varphi, l})$, and $(\mathcal{H}^{\Omega})$ be satisfied. Then each problem $(P_k)$ has an optimal solution provided that $k\in\N$ is sufficiently large.
\end{proposition}\vspace*{-0.1in}
{\bf Proof}. First we show that wherever $k$ is large enough, the discrete trajectory $\{x^{k}_{j}\,|\,j=0\ldots, k\} $ constructed in
Theorem~\ref{strong approximation} for the given local minimizer $\ox(\cdot)$ is a feasible solution to $ (P_{k}) $, i.e., it satisfies all the constraints in \eqref{e:5.1}--\eqref{e:5.3}. Indeed, for the case of \eqref{e:5.1}, this follows directly from \eqref{e:xk2}.
Taking $ k $ such that $ \zeta_{k}\leq\dfrac{\varepsilon}{2} $, we also deduce \eqref{e:5.2} from
\eqref{e:xk2}. By \eqref{beta}, the constraint is \eqref{e:5.3} can be written in the form
\begin{equation*}
\begin{aligned}
\sum\limits_{j=0}^{k-1}\int\limits_{t^{k}_{j}}^{t^{k}_{j+1}}\Big\lVert \dfrac{x^{k}_{j+1}-x^{k}_{j}}{h_{k}}-\dot{\bar{x}}(t)\Big\rVert^{2}\,\mathrm{d}t=\int_{0}^{T}\lVert \dot{x}^{k}(t)-\dot{\bar{x}}(t) \rVert^{2}\,\mathrm{d}t\leq\beta_{k}\leq \frac{\varepsilon}{2},
\end{aligned}
\end{equation*}
and hence holds for large $k$ by the definition of $\beta_{k} $ in \eqref{beta}. Thus $ x^{k}(\cdot) $ is a feasible solution to $ (P_{k}) $ for all large $k$, which tells us that this set is nonempty and also bounded by \eqref{e:5.2}. Now we fix $k\in\N$ and show
that set of feasible solutions to $(P_k)$ is closed. To do this, pick a sequence $z_{\nu}=(x_{\nu}^{0},\ldots,x_{\nu}^{k})$ of feasible solutions to $(P_{k})$ converging to some $z=(x^{0},\ldots,x^{k})$ as $ \nu\to 0 $ and check that $z $ is feasible to $ (P_{k})$.
Indeed, it follows from the closed values of $ F $, the continuity of $ g $, and  the dominated convergence theorem that all the conditions in \eqref{e:5.0}--\eqref{e:5.3}  are satisfied for $z$, which yields the existence of optimal solutions to $(P_k)$ by the Weierstrass theorem in finite dimensions. $\h$

The next theorem establishes a crucial result for the method of discrete approximations in integro-differential systems showing that {\em any} sequences of piecewise-linear extensions of optimal solutions to $(P_k)$ strongly converges in $W^{1,2}$ topology on $[0,T]$
to the prescribed r.i.l.m.\ for $(P)$. Besides being of its own interest, the theorem below opens the door to derive necessary optimality conditions in $(P)$ by passing to the limit from those for its discrete approximations $(P_k)$.\vspace*{-0.1in}

\begin{theorem}[strong convergence of extended discrete optimal
solutions]\label{disc:conver} Let $ \bar{x}(\cdot) $ be an r.i.l.m.
for problem $(P)$. In addition to the standing assumptions
$(\mathcal{H}^{F})$, $(\mathcal{H}^{g})$, $(\mathcal{H}^{\varphi,
l})$, and $(\mathcal{H}^{\Omega})$ imposed along $\bar{x}(\cdot)$,
suppose that the cost functions $ \varphi $ and $ l $ are continuous
at $ \bar{x}(T) $ and at $(t,\bar{x}(\cdot),\dot{\bar{x}}(\cdot)) $
for a.e. $ t\in[0,T] $, respectively. Then any sequence of optimal
solutions $ \bar{x}^{k}(\cdot) $ to $ (P_{k}) $, piecewise linearly
extended on $[0,T] $, converges to $  \bar{x}(\cdot) $ strongly in $
W^{1,2}([0,T],\mathbb{R}^{n}) $ as $k\to\infty$.
\end{theorem}\vspace*{-0.1in} {\bf Proof}. Take any optimal
solutions $\ox^k(\cdot)$ to problem $(P_k)$, which exists for large
$k$ by Proposition~\ref{discri:opti}, and extend it piecewise
linearly to $[0,T]$. We intend to show that
\begin{equation}\label{e:5.4}
\lim\limits_{k\to\infty}\int_{0}^{T}\lVert\dot{\bar{x}}^{k}(t)-\dot{\bar{x}}(t)\rVert^{2}\,\mathrm{d}t=0,
\end{equation} which clearly yields the convergence of $
\bar{x}^{k}(\cdot)$ to $ \bar{x}(\cdot) $ as $k\to\infty$ in the norm topology of $
W^{1,2}([0,T],\mathbb{R}^{n})$. To verify
\eqref{e:5.4}, suppose on the contrary that it fails, i.e., the
limit in \eqref{e:5.4}, along a subsequence (without relabeling)
equals to some $\xi>0$. Indeed, the weak compactness of the unit
ball in $L^{2}([0,T],\mathbb{R}^{n})$ allows us to find $
\upsilon^{x}(\cdot) $ such that $ \dot{\bar{x}}^{k}(\cdot) $
converges weakly to  $ \upsilon^{x}(\cdot)\in
L^{2}([0,T],\mathbb{R}^{n}) $ as $k\to\infty$. According to Mazur's
weak closure theorem, there exists a sequence of convex combinations
of
 $ \dot{\bar{x}}^{k}(\cdot) $ that converges to $ \upsilon^{x}(\cdot)$ in the norm topology of $L^{2}([0,T],\mathbb{R}^{n})$. Hence it contains a subsequence converging to $ \upsilon^{x}(\cdot) $ for a.e. $ t\in [0,T] $. Define $ \tilde{x}(\cdot) $ by
\begin{equation*} \tilde{x}(t):=x_{0}+\int_{0}^{t}
\upsilon^{x}(s)\,\mathrm{d}s,\quad t\in[0,T], \end{equation*} and
get that $ \dot{\tilde{x}}(t)=\upsilon^{x}(t) $ for a.e.
$t\in[0,T]$. By passing to the limit in the conditions of
Theorem~\ref{strong approximation} with the usage of the definition
of $(D_{\mathrm{co}\,F,g}) $ in \eqref{co:inclusion} and the
dominated convergence theorem shows that the limiting  arc
$\tilde{x}(\cdot)$ satisfies the convexified integro-differential
inclusion $ (D_{\mathrm{co}\,F,g}) $, and hence this trajectory is
feasible to the relaxed problem $(R)$.\vspace*{-0.05in}

Next we check that $\tilde{x}(\cdot)$ satisfies the localization conditions in \eqref{loc} relative to $\bar{x}(\cdot)$. Indeed, the first condition in \eqref{loc} follows directly from the passage to the limit in \eqref{e:5.2} as $k\to\infty$. To justify the second
condition in \eqref{loc}, we pass to the limit in \eqref{e:5.3} due to the established weak convergence of the derivatives $\dot{\bar{x}}^{k}(\cdot)\to \dot{\tilde{x}}(\cdot)$ and the lower semicontinuity of the  integral functional\\
$ I[v]:=\int_{0}^{T}\Vert v(t)-\dot{\bar{x}}(t)
\Vert^{2}\,\mathrm{d}t $ in in the weak topology of $L^2$. This
tells us that \begin{equation*} \begin{aligned}
\int_{0}^{T}\lVert\dot{\tilde{x}}(t)-\dot{\bar{x}}(t)\rVert^{2}\leq\liminf\limits_{k\to
\infty}\int_{0}^{T}\Vert\dot{\bar{x}}^{k}(t)-\dot{\tilde{x}}(t)\Vert^{2}\,\mathrm{d}t=\liminf\limits_{k\to
\infty}\sum\limits_{j=0}^{k-1}\int\limits_{t^{k}_{j}}^{t^{k}_{j+1}}\left\lVert
\dfrac{\bar{x}^{k}_{j+1}-\bar{x}^{k}_{j}}{h_{k}}-\dot{\bar{x}}(t)\right\rVert^{2}\,\mathrm{d}t\leq\dfrac{\varepsilon}{2},
\end{aligned} \end{equation*} i.e., we get \eqref{loc}. Moreover, it
follows from the construction of the relaxed problem $(R)$ in
\eqref{R} by the convexity of $\Hat{l}_G$ in the velocity variables,
the weak convergence of the extended discrete derivatives, and the
application  of Mazur's theorem that
\begin{equation*}
\begin{aligned}
\int_{0}^{T}\hat{l}_{G}\big(t,\tilde{x}(t),\dot{\tilde{x}}(t)\big)\,\mathrm{d}t\leq \liminf\limits_{k\to \infty}h_{k}\sum\limits_{j=0}^{k-1}l\Big(t^{k}_{j},\bar{x}^{k}_{j},\dfrac{\bar{x}^{k}_{j+1}
-\bar{x}^{k}_{j}}{h_{k}}\Big).
\end{aligned}
\end{equation*}
We also observe by taking into account
the above definition of $\xi$ that
\begin{equation}\label{e:5.5}
\begin{aligned}
\hat{J}[\tilde{x}]+\dfrac{\xi}{2}=\varphi\big(\tilde{x}(T)\big)+\int_{0}^{T}\hat{l}_{G}\big(t,\tilde{x}(t),\dot{\tilde{x}}(t)\big)\,dt+\dfrac{\xi}{2}\leq\liminf\limits_{k\to
\infty} J_{k}(\bar{x}^{k}). \end{aligned}
\end{equation}
Employing
now Theorem~\ref{strong approximation} produces a sequence
$\{x^{k}(\cdot)\}$ of feasible solutions to $(P_{k})$  that strongly
approximate $\bar{x}(\cdot)$ in $W^{1,2}$, which is a feasible
solution to $(P)$. The forms of the cost functionals in
\eqref{Bolza-functional} and \eqref{disc-bolza} together with the
imposed continuity assumptions on $ \varphi $ and $ l $ yield
\begin{equation}\label{e:5.6} \lim\limits_{k\to \infty}
J_{k}[x^{k}]=J[\bar{x}]. \end{equation} On the other hand, it
follows from the optimality of $\bar{x}^{k}(\cdot)$ in $(P_k)$ that
\begin{equation}\label{e:5.7} J_{k}[\bar{x}^{k}]\leq
J_k[x^{k}]\;\mbox{ for each }\;k\in\mathbb{N}.
\end{equation}
Combining the relationships in \eqref{e:5.5}--\eqref{e:5.7} tells us
that \begin{equation*}
\hat{J}[\tilde{x}]<\hat{J}[\tilde{x}]+\dfrac{\xi}{2}\leq
J[\bar{x}]=\hat{J}[\bar{x}]. \end{equation*} By the above choice of
$\xi>0$, the latter contradicts the fact that $ \bar{x}(\cdot) $ is
an r.i.l.m. for problem $(P)$. Therefore, we get $ \xi=0 $, which
justifies \eqref{e:5.4} and hence completes the proof of the
theorem. $\h$ \vspace*{-0.2in}

\section{Necessary Optimality Conditions for Discrete
Problems}\label{sec:disc-nec}
\setcounter{equation}{0}\vspace*{-0.1in}

According to our scheme of deriving necessary optimality conditions
for local minimizers in the original integro-differential problem
$(P)$, we establish in this section necessary conditions for optimal
solutions to each discrete-time problem $(P_k)$, $k\in\N$, defined
in \eqref{disc-bolza}--\eqref{e:5.3}. Such conditions, being of
their own importance to solve $(P_k)$, can be viewed---by
Theorem~\ref{disc:conver} when $k$ is sufficiently large---as {\em
suboptimal conditions} for the original problem $(P)$, which may be
sufficient for developing numerical methods and practical
applications. Nevertheless, we employ them in the next section to
derive necessary optimality conditions for local minimizers of $(P)$
by passing to the limit as $k\to\infty$.

Fix $k\in\N$ and observe that due to the discrete time in this dynamic
problem, it can be reduced to the following nondynamic problem of
{\em mathematical programming}:
 \begin{equation}\label{mp}
\mathrm{(MP)}\,:\;\left\{ \begin{array}{l}
\text{minimize}\quad \phi_{0}(z)\quad \text{subject to}\vspace{0.2cm}\\
\phi_{j}(z)\leq 0,\quad j=1, \ldots ,s,\vspace{0.2cm}\\
e_{j}(z)=0,\quad j=0,\ldots, m,\vspace{0.2cm}\\
z\in E_{j}\subset\mathbb{R}^{d},\quad j=0,\ldots, l,
\end{array}\right.
\end{equation} where $ \phi_{j}
:\mathbb{R}^{d}\rightarrow\oR $ for $ j=0,\ldots ,s $ and $ e_{j} :
\mathbb{R}^{d}\rightarrow\mathbb{R}^{n} $ for $ j=0,\ldots ,m $. It
is important to emphasize that problem $(MP)$, as well as the
discrete-time problem $(P_k)$, are {\em intrinsically nonsmooth}
even if the functions  $\phi_j$ and $e_j$ are continuously
differentiable, as well as all the functions in $(P)$ and $P_k)$.
The unavoidable nonsmoothness in \eqref{mp} comes from the
increasingly many geometric constraint generated by the {\em
graphical sets} \eqref{e:5.0} in $(P_k)$ that are discretizations of
the velocity sets in the original model \eqref{inclusion}. We know
from Section~\ref{sec:2va} that such constraints are highly
challenging from the viewpoint of generalized differentiation and
cannot be handled by convexified constructions. The robust limiting
constructions with full calculus overviewed in Section~\ref{sec:2va}
are the most appropriate for our purposes.

The following necessary optimality conditions for problem (MP) are
consequences of \cite[Theorem~6.5(iii)]{b3} and the intersection
rule for limiting normals taken from
\cite[Corollary~2.17]{b3}.\vspace*{-0.1in}

\begin{lemma}[\bf necessary optimality conditions for mathematical
programs]\label{p:mp} Let $ \bar{z} $ be an optimal solution to
problem \eqref{mp}. Assume also that $ \phi_{j}$ for $j=0,1,\ldots,s$
are locally Lipschitzian around $ \bar{z} $, that $ e_{j}$ for
$j=1,\ldots,l$ are ${\cal C}^1$-smooth around this point, and that
each set $E_{j}$ for $j=1,\ldots,m$ is closed around $ \bar{z} $.
Then there are real numbers $
\{\mu_{j}\in\mathbb{R}\;|\;j=0,\ldots,s\} $ as well as vectors $
\{z^{*}_{j}\in \mathbb{R}^{d}\;|\;j=0,\ldots,l\} $ and $
\{p_{j}\in\mathbb{R}^{n}\;|\;j=0,\ldots,m\} $, not equal to zero
simultaneously, such that $ \mu_{j}\geq 0 $ for $ j=0,\ldots,s $,
and we have \begin{equation}\label{muj}
\mu_{j}\phi_{j}(\bar{z})=0\quad\text{for}\;\;j=1,\ldots,s,
\end{equation} \begin{equation*}z^{*}_{j}\in
N_{E_j}(\bar{z})\quad\text{for}\;\;j=0,\ldots,l, \end{equation*}
\begin{equation*}
 -\sum\limits_{j=0}^{l}z_{j}^{*}\in
\partial\,\Big(\sum\limits_{j=0}^{s}\mu_{j}\phi_{j}\Big)(\bar{z})+\sum\limits_{j=0}^{m}\Big(\nabla\,e_{j}(\bar{z})\Big)^{*}p_{j}.
\end{equation*}
 \end{lemma}

Now we employ Lemma~\ref{p:mp} and calculus rules for the
generalized differential constructions from Section~\ref{sec:2va} to
derive necessary optimality conditions in problems $ (P_{k}) $ for
each fixed $ k\in\N $.\vspace*{-0.1in}

\begin{theorem}[\bf necessary conditions in discrete
problems]\label{discr:N.C} Fix $k\in\N$ to be sufficiently large,
and let $ (x_{0},\ldots,\bar{x}^{k}_{k}) $ be an optimal solution to
problem $(P_k)$. Suppose that the cost functions $\ph$ and $l$  are
locally Lipschitzian around the corresponding components of the
optimal solution for all $t\in\Delta_k$ $ ($with the index $t$
dropped when no confusion arise$)$, that the function $g(t,s,\cdot)$
is continuously differentiable around the optimal points, and that
the multifunction $ F(t_{j},\cdot) $ is bounded and Lipschitz
continuous around $ \bar{x}^{k}_{j} $ with the sets ${\rm
gph}\,F(t_{j},\cdot) $ and $ \Omega $ being locally closed. Then
there exist a number $ \lambda^{k}\geq 0 $ and a vector $
p^{k}=(p^{k}_{0},\ldots,p^{k}_{k})\in\mathbb{R}^{(k+1)n} $ such that
\begin{equation}\label{discr:eq1} \lambda^{k}+\Vert p^{k}_{k}\Vert
\neq 0, \end{equation} \begin{equation}\label{discr:eq2}
\begin{aligned}
-p^k_{k}\in\lambda^{k}\partial\varphi\big(\bar{x}^{k}_{k})+
N_{\Omega_{k}}(\bar{x}^{k}_{k}),
\end{aligned}
\end{equation}
\begin{equation}\label{discr:eq3} \begin{aligned}
&\Big(\dfrac{p^k_{j+1}-p^k_{j}}{h_{k}}+2h_{k}^{-1}\mu_{j}p^k_{j+1}-h_{k}^{-1}\mu_{j}\lambda^{k}(v_{j}+h_{k}^{-1}\theta_{j})+h_{k}^{-1}\sum\limits_{\nu=j+1}^{k}\tilde{\xi}^{j}_{\nu}p^k_{\nu+1} , p^k_{j+1}-\lambda^{k}h_{k}^{-1}\theta_{j}\Big)\\
&\in\lambda^{k}\partial\,l\Big(t_{j},\bar{x}_{j},\dfrac{\bar{x}_{j+1}-\bar{x}_{j}}{h_{k}}\Big)+
N_{{\rm\small
gph}\,F_{j}}\Big(\bar{x}^{k}_{j},\frac{\bar{x}^{k}_{j+1}-\bar{x}^{k}_{j}}{h_{k}}-\bar{w}^{k}_{j}\Big),\;\;j=0,\ldots,k-1,
\end{aligned}
\end{equation}
where $ F_{j}(\cdot):=F(t_{j},\cdot) $,
and where we use the notations
\begin{equation}\label{eqq:wkj}
\bar{w}^{k}_{j}=\dfrac{1}{h^{k}_{j}}\int\limits_{t^{k}_{j}}^{t^{k}_{j+1}}\Big\{\sum\limits_{i=0}^{j-1}\int\limits_{t^{k}_{i}}^{t^{k}_{i+1}}g(t,s,\bar{x}^{k}_{i})\,\mathrm{d}s+\int\limits_{t^{k}_{j}}^{t}g(t,s,\bar{x}^{k}_{j})\,\mathrm{d}s\Big\}\,\mathrm{d}t,
\;\;j=0,\ldots,k-1,
\end{equation}
\begin{equation}\label{theta}
\begin{aligned} \theta_{j}:= \int\limits_{t_{j}}^{t_{j+1}}\Big(
\dfrac{\bar{x}_{j+1}-\bar{x}_{j}}{h_{k}}-\dot{\bar{x}}(t)\Big)\,\mathrm{d}t,\;\;j=0,\ldots,k-1,
\end{aligned} \end{equation} \begin{equation}\label{xi} \begin{aligned}
\xi^{j}_{i}:=\int_{t_{i}}^{t_{i+1}}\Big\{\int_{t_{j}}^{t_{j+1}}\nabla\,g(t,s,\bar{x}_{j})^*\,\mathrm{d}s\Big\}\mathrm{d}t,\;\;0\leq
j\leq i-1,\;i=1,\ldots,k-1, \end{aligned} \end{equation}
\begin{equation}\label{e:xii}
\tilde{\xi}_{i}^{j}:=\left\{\begin{array}{lll}
\xi_{i}^{j} & \mbox{if}\quad 1\leq i\leq k-1\\
0 &\mbox{if}\quad i=k \end{array}\right.,\quad 0\leq j\leq
i-1,\;\;1\leq i\leq k, \end{equation}
\begin{equation}\label{mu}
\begin{aligned}
\mu_{j}:=\int_{t_{j}}^{t_{j+1}}\Big\{\int_{t_{j}}^{t}\nabla\,g(t,s,\bar{x}_{j})^*\,\mathrm{d}s\Big\}\mathrm{d}t,\;\;j=0,\ldots,k-1,
\end{aligned} \end{equation}
with $\nabla g(t,s,\cdot)^*$ standing for the adjoint Jacobian of $g$ with respect to the state variable.
\end{theorem}\vspace*{-0.1in} {\bf
Proof}. It what follows, we write $ t^{k}_{j}:=t_{j} $ (resp. $
\upsilon^{k}_{j}:=\upsilon_{j} $) for the mesh points (resp.
approximations) as $ j=0,\ldots,k $. Denote
$z:=(x_{0},\ldots,x_{k},X_{0},\ldots,X_{k-1}) $, where $ x_{0} $ is
fixed, and consider the mathematical program  of type \eqref{mp}
given by
\begin{equation}\label{MP:1}
\begin{aligned}
\text{minimize}\quad\phi_{0}(z):=\varphi(x^{k}_{k})+h_{k}\sum\limits_{j=0}^{k-1}l(t_{j},x_{j},X_{j})+\dfrac{1}{2}\sum\limits_{j=0}^{k-1}\int\limits_{t_{j}}^{t_{j+1}}\lVert
X_{j}-\dot{\bar{x}}(t)\rVert^{2}\,\mathrm{d}t,
\end{aligned}
\end{equation} subject to the inequality, equality, and geometric
constraints \begin{equation*}\label{MP:2} \phi_{j}(z):=\lVert
x_{j-1}-\bar{x}(t_{j-1}) \rVert-\dfrac{\varepsilon}{2}\leq 0,\quad
j=1,\ldots,k+1, \end{equation*} \begin{equation*}\label{MP:3}
\phi_{k+2}(z):=\sum\limits_{j=0}^{k-1}\int\limits_{t_{j}}^{t_{j+1}}\lVert
X_{j}-\dot{\bar{x}}(t)\,\mathrm{d}t-\dfrac{\varepsilon}{2}\leq 0,
\end{equation*} \begin{equation*}\label{MP:4}
e_{j}(z):=x_{j+1}-x_{j}-\int\limits_{t_{j}}^{t_{j+1}}\Big\{\sum\limits_{i=0}^{j-1}\int\limits_{t_{i}}^{t_{i+1}}g(t,s,x_{i})\,\mathrm{d}s+\int\limits_{t_{j}}^{t}g(t,s,x_{j})\,\mathrm{d}s\Big\}\mathrm{d}t-h_{k}X_{j}=0,\quad
j=0,\ldots,k-1, \end{equation*} \begin{equation}\label{MP:6} z\in
E_{j}:=\big\{z\in\mathbb{R}^{(2k+1)n}\;\big|\; X_{j}\in
F(t_{j},x_{j})+w_{j}\big\},\quad j=0,\ldots,k-1, \end{equation}
\begin{equation}\label{MP:7} z\in
E_{k}:=\big\{z\in\mathbb{R}^{(2k+1)n}\;\big|\;x_{0}\;\;\mbox{is
fixed},\;\;x^{k}_{k}\in \Omega_{k}\big\}. \end{equation}
Apply now
the necessary optimality conditions from Lemma~\ref{p:mp} with $
s=k+2 $, $ m=k-1 $ and $ l=k $ to the optimal solution $
\bar{z}^{k}=(\bar{x}^{k}_{0},\ldots,\bar{x}^{k}_{k},\bar{X}^{k}_{0},\ldots,\bar{X}^{k}_{k-1})
$ of \eqref{mp}, where $
\bar{x}^{k}=(\bar{x}^{k}_{0},\ldots,\bar{x}^{k}_{k}) $ is a given
optimal solution of $(P_{k}) $. First observe that the inequality
constraints in \eqref{mp} defined by the functions $\phi_j$ with
$j=1,\ldots,k$ are inactive when $k$ sufficiently large due to the
$W^{1,2}$-strong convergence $\bar{x}^k(\cdot)\to\bar{x}(\cdot)$ as
$k\to\infty$ established above in Theorem~\ref{strong
approximation}, which means that $ \phi_{j}(\bar{z}^{k})<0 $ as $
j=1,\ldots,k+2 $. All of this and the complementary slackness
conditions in \eqref{muj} allow us to find $
\mu_{0}:=\lambda^{k}\geq 0 $ with $ \mu_{j}=0 $ for $ j=1,\ldots,k+2
$, $ p_{j}^{k}\in \mathbb{R}^{n} $ for $j=1,\ldots,k$, and
\begin{equation*} \begin{aligned} z^{\ast}_{j}=\big(&
x^{\ast}_{0j},\ldots,x^{\ast}_{kj},X^{\ast}_{0j},
\ldots,X^{\ast}_{(k-1)j}\big),\;\;j=0,\ldots,k-1, \end{aligned}
\end{equation*} which are not equal to zero simultaneously and such
that
\begin{equation}\label{MP:8} z^{\ast}_{j}\in
N_{E_j}(\bar{z}^{k}),\;\;j=0,\ldots,k, \end{equation}
\begin{equation*}\label{MP:9} -\sum\limits_{j=0}^{k}z^{\ast}_{j}\in
\lambda^{k}\partial\phi_{0}(\bar{z}^{k})+\sum\limits_{j=0}^{k-1}\left(\nabla\,e_{j}(\bar{z}^{k})\right)^{*}p_{j+1},
\end{equation*} where we skip the upper index `$k$' for $p$ to simplify the notation. It follows from the structures of the sets $ E_{j} $
in \eqref{MP:6} and \eqref{MP:7} that the normal cone inclusions in
\eqref{MP:8} can be written in the form \begin{equation}\label{e:Fj}
(x^{*}_{jj},X^{*}_{jj})\in N_{{\rm\small
gph}\,G_{j}}\Big(\bar{x}^{k}_{j},\frac{\bar{x}^{k}_{j+1}-\bar{x}^{k}_{j}}{h_{k}}\Big),\;\;j=0,\ldots,k-1,
\end{equation} where $G_j$ reflects the discrete
integro-differential dynamics by \begin{equation*}
G_{j}(\cdot):=F(t_{j},\cdot)+\dfrac{1}{h^{k}_{j}}\int\limits_{t^{k}_{j}}^{t^{k}_{j+1}}\Big\{\sum\limits_{i=0}^{j-1}\int\limits_{t^{k}_{i}}^{t^{k}_{i+1}}g(t,s,\bar{x}^{k}_{i})\,\mathrm{d}s+\int\limits_{t^{k}_{j}}^{t}g(t,s,\cdot)\,\mathrm{d}s\Big\}\,\mathrm{d}t,
\;\;j=0,\ldots,k-1, \end{equation*} and where the corresponding
components of the dual vectors $z^*_j$ satisfy the conditions in \eqref{e:Fj} and
\begin{equation}\label{xkj=0}
x^{*}_{ij}=X^{*}_{ij}=0\;\;\text{if}\;\;i\neq
j\;\;\text{for}\;\;j=0,\ldots,k-1,
\end{equation}
\begin{equation}\label{x0k-xkk} x^{*}_{kk}\in N_{\Omega_{k}}(\bar{x}^{k}_{k}) \end{equation}
with $x^*_{0,k}$ being free by \eqref{MP:7}. Therefore, we arrive at the representations
\begin{equation*}
\begin{aligned}
&-\sum\limits_{j=0}^{k}z^{\ast}_{j}=(-x^{\ast}_{00}-x^{\ast}_{0k},-x^{\ast}_{11},\ldots,-x^{\ast}_{(k-1)(k-1)},-x^{\ast}_{kk},-X^{\ast}_{00},\ldots,-X^{\ast}_{(k-1)(k-1)}),\\
&\quad\sum\limits_{j=0}^{k-1}\left(\nabla\,e_{j}(\bar{z}^{k})\right)^{*}p_{j+1}\\
&\quad=(-p_{1}-\mu_{0}p_{1}-\sum\limits_{j=1}^{k-1}\xi^{0}_{j}p_{j+1}, p_{1}-p_{2}-\mu_{1}p_{2}-\sum\limits_{j=2}^{k-1}\xi^{1}_{j}p_{j+1}, p_{2}-p_{3}-\mu_{2}p_{3}-\sum\limits_{j=3}^{k-1}\xi^{2}_{j}p_{j+1},\\
&\quad\ldots,
p_{k-2}-p_{k-1}-\mu_{k-2}p_{k-1}-\xi^{k-2}_{k-1}p_{k},
p_{k-1}-p_{k}-\mu_{k-1}p_{k},p_{k},-h_{k}p_{1},-h_{k}p_{2},\ldots,-h_{k}p_{k})
\end{aligned}
\end{equation*} with $\xi^j_k$ and $\mu_j$ defined in \eqref{xi} and \eqref{mu}, respectively.
On the other hand, applying the basic subdifferential sum rule from \cite[Theorem~2.19]{b3} to
the minimizing function in \eqref{MP:1} with taking into account its
structure gives us the inclusion
\begin{equation}\label{e:6.13} \begin{aligned}
\partial \phi_{0}(\bar{z}^{k})\subset\partial
\varphi(\bar{x}^{k}_{k})+h_{k}\sum\limits_{j=0}^{k-1}\partial
l\big(t_{j},\bar{x}_{j},\bar{X}_{j}\big)+\sum\limits_{j=0}^{k-1}\nabla
\rho_{j}(\bar{z}^{k}) \end{aligned} \end{equation} with $
\rho_{j}(\bar{z}^{k}):=\dfrac{1}{2}\displaystyle\int\limits_{t_{j}}^{t_{j+1}}\lVert
\bar{X}_{j}-\dot{\bar{x}}(t)) \rVert^{2}\,\mathrm{d}t$ and
with\vspace*{-0.1in} \begin{enumerate} \item[$  \bullet $] $
\sum\limits_{j=0}^{k-1}\nabla
\rho_{j}(\bar{z}^{k})=(0,\ldots,0,\theta_{0},\theta_{1},\ldots,\theta_{k-1}),
$ where $\theta_{j}$ is defined in \eqref{theta}, \item[$ \bullet $]
$ \sum\limits_{j=0}^{k-1}\partial
l(t_{j},\bar{x}_{j},\bar{X}_{j})=(w_{0},\ldots,w_{k},v_{0},\ldots,v_{k-1})$,
where $ (w_{j},v_{j})\in
\partial\,l\Big(t_{j},\bar{x}_{j},\dfrac{\bar{x}_{j+1}-\bar{x}_{j}}{h_{k}}\Big)$,
\item[$ \bullet $]
$\partial\varphi(\bar{x}^{k}_{k})=(0,\ldots,u^{k},\ldots,0) $.
\end{enumerate} This allows us to represent the term $
\lambda^{k}\partial\phi_{0}(\bar{z}^{k}) $ in \eqref{e:6.13} as
\begin{equation*}
\lambda^{k}(h_{k}w_{0},\ldots,h_{k}w_{k-1},u^{k},\theta_{0}+h_{k}v_{0},\ldots,\theta_{k-1}+h_{k}v_{k-1}).
\end{equation*} Combining all the above decomposes the inclusion in
\eqref{e:6.13} into the following equalities:
\begin{equation}\label{e:6.19}
-x^{\ast}_{00}-x^{\ast}_{0k}=\lambda^{k}h_{k}w_{0}-p_{1}-\mu_{0}p_{1}-\sum\limits_{j=1}^{k-1}\xi^{0}_{j}p_{j+1},
\end{equation} \begin{equation}\label{e:6.16}
-x^{\ast}_{jj}=\lambda^{k}h_{k}w_{j}+p_{j}-p_{j+1}-\mu_{j}p_{j+1}-\sum\limits_{\nu=j+1}^{k-1}\xi^{j}_{\nu}p_{\nu+1},\,\,\,j=1,\ldots,k-2,
\end{equation} \begin{equation}\label{e:6.20}
-x^{\ast}_{(k-1)(k-1)}=\lambda^{k}h_{k}w_{k-1}+p_{k-1}-p_{k}-\mu_{k-1}p_{k},
\end{equation} \begin{equation}\label{e:6.17}
-x^{\ast}_{kk}=\lambda^{k}u^{k}+p_{k}, \end{equation}
\begin{equation}\label{e:6.18}
-X^{\ast}_{jj}=\lambda^{k}h_{k}v_{j}+\lambda^{k}\theta_{j}-h_{k}p_{j+1},\,\,\,j=0,\ldots,k-1.
\end{equation} Denoting $ p_{0}:=x^{*}_{0k} $, $p_{k+1}:=0 $ and
taking $\tilde{\xi}_{i}^{j}$ in \eqref{e:xii}, we deduce from
\eqref{e:6.19}, \eqref{e:6.16}, and \eqref{e:6.18} that
\begin{equation*} \begin{aligned}
\Big(\dfrac{p_{j+1}-p_{j}}{h_{k}}-\lambda^{k}w_{j}+h_{k}^{-1}\mu_{j}p_{j+1}+h_{k}^{-1}\sum\limits_{\nu=j+1}^{k}\tilde{\xi}^{j}_{\nu}p_{\nu+1},
p_{j+1}-\lambda^{k}(v_{j}+h_{k}^{-1}\theta_{j})\Big)=h_{k}^{-1}(x^{*}_{jj}).
\end{aligned} \end{equation*} Moreover, employing \eqref{e:Fj} tells us
that\begin{equation*} \begin{aligned}
&\Big(\dfrac{p_{j+1}-p_{j}}{h_{k}}-\lambda^{k}w_{j}+h_{k}^{-1}\mu_{j}p_{j+1}+h_{k}^{-1}\sum\limits_{\nu=j+1}^{k}\tilde{\xi}^{j}_{\nu}p_{\nu+1} , p_{j+1}-\lambda^{k}(v_{j}+h_{k}^{-1}\theta_{j})\Big)\\
&\in N_{{\rm\small
gph}\,G_{j}}\Big(\bar{x}^{k}_{j},\frac{\bar{x}^{k}_{j+1}-\bar{x}^{k}_{j}}{h_{k}}\Big),\quad
j=0,\ldots,k-1, \end{aligned} \end{equation*} which is equivalent,
by the coderivative definition in \eqref{cod}, to \begin{equation*}
\begin{aligned}
&\dfrac{p_{j+1}-p_{j}}{h_{k}}-\lambda^{k}w_{j}+h_{k}^{-1}\mu_{j}p_{j+1}+h_{k}^{-1}\sum\limits_{\nu=j+1}^{k}\tilde{\xi}^{j}_{\nu}p_{\nu+1} \\
&\in
D^{*}G_{j}\big(\bar{x}^{k}_{j},(\bar{x}^{k}_{j+1}-\bar{x}^{k}_{j})/h_{k}\big)\big(\lambda^{k}(v_{j}+h_{k}^{-1}\theta_{j})-p_{j+1}\big),\quad
j=0,\ldots,k-1. \end{aligned} \end{equation*} Due to the summation
structure of $G_j$ and the Lipschitz continuity of $F_j$, we apply
to $G_j$ the coderivative sum rule from \cite[Theorem~3.9]{b3}
together with Leibniz's rule of differentiation under the integral
sign while getting in this way that \begin{equation*} \begin{aligned}
&\dfrac{p_{j+1}-p_{j}}{h_{k}}-\lambda^{k}w_{j}+h_{k}^{-1}\mu_{j}p_{j+1}+h_{k}^{-1}\sum\limits_{\nu=j+1}^{k}\tilde{\xi}^{j}_{\nu}p_{\nu+1}\\&\in
D^{*}F_{j}\big(\bar{x}^{k}_{j},(\bar{x}^{k}_{j+1}-\bar{x}^{k}_{j})/h_{k}-\bar{w}^{k}_{j}\big)\big(\lambda^{k}(v_{j}
+h_{k}^{-1}\theta_{j})-p_{j+1}\big)\\
&+\dfrac{1}{h^{k}_{j}}\Big[\int\limits_{t^{k}_{j}}^{t^{k}_{j+1}}\Big\{\int\limits_{t^{k}_{j}}^{t}\nabla\,g(t,s,\bar{x}_{j})^*\,\mathrm{d}s\Big\}\,\mathrm{d}t\Big]\,\big(\lambda^{k}(v_{j}+h_{k}^{-1}\theta_{j})-p_{j+1}\big).
\end{aligned} \end{equation*} This allows us to verify that
\begin{equation}\label{graph1} \begin{aligned}
&\dfrac{p_{j+1}-p_{j}}{h_{k}}-\lambda^{k}w_{j}+h_{k}^{-1}\mu_{j}p_{j+1}+h_{k}^{-1}\sum\limits_{\nu=j+1}^{k}\tilde{\xi}^{j}_{\nu}p_{\nu+1}\\&\in
D^{*}F_{j}\left(\bar{x}^{k}_{j},(\bar{x}^{k}_{j+1}-\bar{x}^{k}_{j})/h_{k}-\bar{w}^{k}_{j}\right)(\lambda^{k}(v_{j}+h_{k}^{-1}\theta_{j})-p_{j+1})\\&+h_{k}^{-1}\mu_{j}(\lambda^{k}(v_{j}+h_{k}^{-1}\theta_{j})-p_{j+1}),
\end{aligned} \end{equation} which implies in turn the inclusion
\begin{equation*} \begin{aligned}
&\Big(\dfrac{p_{j+1}-p_{j}}{h_{k}}-\lambda^{k}w_{j}+2h_{k}^{-1}\mu_{j}p_{j+1}-h_{k}^{-1}\mu_{j}\lambda^{k}(v_{j}+h_{k}^{-1}\theta_{j})+h_{k}^{-1}\sum\limits_{\nu=j+1}^{k}\tilde{\xi}^{j}_{\nu}p_{\nu+1} , p_{j+1}-\lambda^{k}(v_{j}+h_{k}^{-1}\theta_{j})\Big)\\
&\in N_{{\rm\small
gph}\,F_{j}}\Big(\bar{x}^{k}_{j},\frac{\bar{x}^{k}_{j+1}-\bar{x}^{k}_{j}}{h_{k}}-\bar{w}^{k}_{j}\Big).
\end{aligned} \end{equation*} The latter clearly yields the
relationship \begin{equation*}\label{graph} \begin{aligned}
&\Big(\dfrac{p_{j+1}-p_{j}}{h_{k}}+2h_{k}^{-1}\mu_{j}p_{j+1}-h_{k}^{-1}\mu_{j}\lambda^{k}(v_{j}+h_{k}^{-1}\theta_{j})+h_{k}^{-1}\sum\limits_{\nu=j+1}^{k}\tilde{\xi}^{j}_{\nu}p_{\nu+1} , p_{j+1}-\lambda^{k}h_{k}^{-1}\theta_{j}\Big)\\
&\in\lambda^{k}\partial\,l\Big(t_{j},\bar{x}_{j},\dfrac{\bar{x}_{j+1}-\bar{x}_{j}}{h_{k}}\Big)+ N_{{\rm\small gph}\,F_{j}}\Big(\bar{x}^{k}_{j},\frac{\bar{x}^{k}_{j+1}-\bar{x}^{k}_{j}}{h_{k}}-\bar{w}^{k}_{j}\Big).
\end{aligned}
\end{equation*}
Employing further \eqref{x0k-xkk} and \eqref{e:6.17}, we get
\begin{equation*}
\begin{aligned}
-p_{k}=\lambda^{k}u^{k}+ x^{*}_{kk}\in\lambda^{k}\partial\varphi(\bar{x}^{k}_{k})+N_{\Omega_{k}}(\bar{x}^{k}_{k})
.
\end{aligned}
\end{equation*}
To complete the proof of the theorem, it remains to show that the nontriviality condition \eqref{discr:eq1} is satisfied. Supposing on the contrary that $ \lambda^{k}=0 $ and $ p_{k}^k=0 $ ensures by \eqref{graph1} that
\begin{equation*}
\begin{aligned}
\dfrac{p_{j+1}^k-p_{j}^k}{h_{k}}+h_{k}^{-1}\mu_{j}p_{j+1}^k+2h_{k}^{-1}\sum\limits_{\nu=j+1}^{k}\tilde{\xi}^{j}_{\nu}p_{\nu+1}^k\in D^{*}F_{j}\big(\bar{x}^{k}_{j},(\bar{x}^{k}_{j+1}-\bar{x}^{k}_{j})/h_{k}-\bar{w}^{k}_{j}\big)(-p_{j+1}^k).
\end{aligned}
\end{equation*}
This tells us, by the coderivative criterion for Lipschitz continuity in Theorem~\ref{est:cod}, that $ p^k_{k}=0 $ yields $ p^k_{j}=0 $ for all $ j=0,\ldots,k-1 $,  and hence we get the implications
\begin{equation*}
\begin{aligned}
\eqref{e:6.19}-\eqref{e:6.17}\Longrightarrow x_{jj}^{*}=0,\;\;j=0,\ldots,k,\\
\eqref{e:6.18}\Longrightarrow X_{jj}^{*}=0,\;\;j=0,\ldots,k-1.
\end{aligned}
\end{equation*}
We know from \eqref{xkj=0} that all the components of $ z^{\ast}_{j} $ different from $ (x^{\ast}_{jj},X^{\ast}_{jj}) $ are zero for $ j=0,\ldots,k-1 $. Therefore, $ z^{*}_{j}=0 $ for all $ j=0,\ldots,k $, which contradicts the nontriviality condition in problem \eqref{mp} of mathematical programming formulated in Lemma~\ref{p:mp} and thus completes the proof. $ \h $\vspace*{-0.2in}

\section{Necessary Conditions for Integro-Differential Problems}
\label{sec:bolza-nec} \setcounter{equation}{0}\vspace*{-0.1in}

Now we are at the third and final stage of the method of discrete approximations to derive necessary optimality conditions for the constrained integro-differential inclusion problem $(P)$. This stage accumulates the previous results on the convergence of discrete optimal solutions and the necessary optimality conditions for $(P_k)$ to derive necessary conditions for the prescribed r.i.l.m. $\ox(\cdot)$ of $(P)$ by passing to the limit as $k\to\infty$. The major part of our analysis in this section is to establish an appropriate convergence of {\em dual functions} (adjoint arcs) for which the Lipschitz continuity of the velocity multifunction $F$ plays a crucial role.

To furnish the passage to the limit from the necessary optimality conditions for discrete problems, we exploit the {\em robustness} of our basic generalized differential constructions reviewed in Section~\ref{sec:2va} with respect to perturbations of the initial data. While such a robustness is inherent in our basic constructions relative to their variables, the limiting procedures in what follows include also the time parameter $t$ in the general setting under consideration. To accomplish this, we impose the following rather mild requirements for the case of nonautonomous IDIs, where `Limsup' stands for the (Kuratowski-Painlev\'e) {\em outer limit} of multifunctions; see \cite{b3,rw}:
\begin{enumerate}
\item[$ (\mathcal{A}^{l}) $]\,Given $l\colon[0,T]\times\R^n\times\R^n\to\oR$ and $(\ox,\bar v,\bar t)\in\dom l$, we have
\begin{equation*}
\Limsup\limits_{(t,x,v)\rightarrow(\bar t,\bar{x},\bar{v})}\partial\,l(t,x,v)=\partial\,l(\bar t,\bar{x},\bar{v}),
\end{equation*}
where the subdifferential is taken with respect to $(x,v)$.

\item[$ (\mathcal{A}^{F}) $]\,Given $F\colon[0,T]\times\R^n\tto\R^n$ and $(\ox,\ov)\in\gph F(t,\cdot)$, we have
\begin{equation*}
\limsup\limits_{(x,v,t)\rightarrow(\bar{x},\bar{v},\bar{t})}\,N\big((x,v) ; \mathrm{gph}\,F(\cdot , t)\big)=N\big((\bar{x},\bar{v}); \mathrm{gph}\,F(\cdot , \bar{t})\big).
\end{equation*}
\end{enumerate}
Properties $ (\mathcal{A}^{l}) $ and $ (\mathcal{A}^{F}) $ hold in fairly general settings; in particular, in the separated case when $ l=l_{1}(x,v)+l_{2}(t) $ and  $ F=F_{1}(x)+F_{2}(t) $; see, e.g., \cite[Proposition~5.70]{Boris}.

Here is our principal result about necessary optimality conditions for the Bolza problem governed by nonconvex constrained IDIs. We are not familiar with any previous versions of the obtained {\em extended Volterra condition} type even in particular cases of smooth and unconstrained IDIs.\vspace*{-0.1in}

\begin{theorem}[\bf Volterra-type necessary optimality conditions for constrained IDIs]\label{main-result} Let $\ox(\cdot)$ be a r.i.l.m.\ of the IDI problem $(P)$, where $ (\mathcal{A}^{l}) $ and $ (\mathcal{A}^{F}) $ are satisfied in addition to the assumptions of Theorem~{\rm\ref{discr:N.C}}. Then there exist a multiplier $\lambda\geq 0$ and an absolutely continuous mapping $ p : [0,T]\rightarrow\mathbb{R}^{n} $ such that the following conditions hold:\\[1ex]
{\sc Generalized Volterra Condition}: for a.e. $ \tau\in[0,T] $ we have
\begin{equation}\label{vol-cone}
\begin{aligned}
&\dot{p}(\tau)+\int_{\tau}^{T}\nabla\,g\big(t,\tau,\bar{x}(\tau)\big)^*\,p(t)\,\mathrm{d}t\\
 &\in\mathrm{co}\,\Big\{u\;\Big|\;\big(u , p(\tau)\big)\in\lambda\,\partial\,l\big(\tau , \bar{x}(\tau) , \dot{\bar{x}}(\tau)\big)+ N\big((\bar{x}(\tau),\dot{\bar{x}}(\tau)-\int_{0}^{\tau}g\big(\tau,t,\bar{x}(t)\big)\,\mathrm{d}t);\mathrm{gph}\,F(\tau,\cdot)\big)\Big\},
\end{aligned}
\end{equation}
which can be equivalently written in the coderivative form: for a.e. $ \tau\in[0,T] $ it holds that
\begin{equation}\label{vol-cod}
\begin{aligned}
&\dot{p}(\tau)+\int_{\tau}^{T}\nabla\,g\big(t,\tau,\bar{x}(\tau)\big)^*\,p(t)\,\mathrm{d}t\\
& \in\mathrm{co}\,\left\{\underset{(v,w)\in\partial\,l\big(\tau , \bar{x}(\tau) , \dot{\bar{x}}(\tau)\big)}{\bigcup}\Big[\lambda\,v+D^{\ast}F\Big(\tau,\bar{x}(\tau),\dot{\bar{x}}(\tau)-\int_{0}^{\tau}g\big(\tau,t,\bar{x}(t)\big)\,\mathrm{d}t\Big)\big(\lambda\,w-p(\tau)\big)\Big]\right\}.
\end{aligned}
\end{equation}
{\sc Transversality condition}:
\begin{equation}\label{tran}
\begin{aligned}
-p(T)\in\lambda\,\partial\varphi\big(\bar{x}(T)\big)+N_{\Omega}\big(\bar{x}(T)\big).
\end{aligned}
\end{equation}
{\sc Nontriviality condition}:
\begin{equation}\label{nontriv}
\lambda+\Vert p(T)\Vert =1.
\end{equation}
\end{theorem}\vspace*{-0.1in}
{\bf Proof}. We split the proof into the four major steps as follows.\\[1ex]
{\bf Step~1}. {\em Estimates the approximating dual elements.}
According to \eqref{discr:eq3} and the coderivative definition \eqref{cod}, there exist vectors $(v^{k}_{j},w^{k}_{j})\in \partial\,l\left(t^{k}_{j},\bar{x}^{k}_{j},(\bar{x}^{k}_{j+1}-\bar{x}^{k}_{j})/h_{k}\right) $ such that
\begin{equation*}
\begin{aligned}
&\dfrac{p_{j+1}-p_{j}}{h_{k}}-\lambda^{k}w^{k}_{j}+2h_{k}^{-1}\mu_{j}p_{j+1}-h_{k}^{-1}\mu_{j}\lambda^{k}(v_{j}+h_{k}^{-1}\theta_{j})+h_{k}^{-1}\sum\limits_{\nu=j+1}^{k}\tilde{\xi}^{j}_{\nu}p_{\nu+1}\\
&\quad\quad\in D^{*}F_{j}\big(\bar{x}^{k}_{j},(\bar{x}^{k}_{j+1}-\bar{x}^{k}_{j})/h_{k}-\bar{w}^{k}_{j}\big)\big(-p_{j+1}+\lambda^{k}(v_{j}+h_{k}^{-1}\theta_{j})\big),
\end{aligned}
\end{equation*}
for all $ j=0,\ldots,k-1 $. From the above relationships, Theorem~\ref{strong approximation}, $ (\mathcal{H}^{F}_{2})$, and Theorem~\ref{est:cod} we get
\begin{equation}\label{e:p_prime}
\begin{aligned}
&\Big\Vert\dfrac{p_{j+1}-p_{j}}{h_{k}}-\lambda^{k}w^{k}_{j}+2h_{k}^{-1}\mu_{j}p_{j+1}-h_{k}^{-1}\mu_{j}\lambda^{k}(v_{j}+h_{k}^{-1}\theta_{j})+h_{k}^{-1}\sum\limits_{\nu=j+1}^{k}\tilde{\xi}^{j}_{\nu}p_{\nu+1}\Big\Vert\\
&\quad\quad \leq l_{F}\Vert-p_{j+1}+\lambda^{k}(v_{j}+h_{k}^{-1}\theta_{j})\Vert,\quad  j=0,\ldots,k-1.
\end{aligned}
\end{equation}
It follows from the subdifferential estimate \eqref{e:sub} that
\begin{equation}\label{e:wv}
\Vert w_{j} \Vert\leq M_{l}\;\mbox{ and }\;\Vert v_{j} \Vert\leq M_{l}\;\;\text{whenever}\;\;j=0,\ldots,k-1.
\end{equation}
Using \eqref{e:p_prime} and \eqref{e:wv}, we arrive at the chain of inequalities
\begin{equation*}
\begin{aligned}
\Vert p_{j} \Vert &\leq (1+h_{k}l_{F}+2\Vert\mu_{j}\Vert)\Vert p_{j+1}\Vert+\sum\limits_{\nu=j+1}^{k}\Vert\tilde{\xi}^{j}_{\nu}\Vert\cdot \Vert p_{\nu+1}\Vert\\
& + \lambda^{k}M_{l}(h_{k}(1+l_{F})+\Vert\mu_{j}\Vert)+\lambda^{k}(h_{k}^{-1}\Vert \mu_{j}\Vert+l_{F})\Vert \theta_{j}\Vert
\\&\leq (1+h_{k}l_{F}+2\alpha h_{k}^{2})\Vert p_{j+1}\Vert+\alpha h_{k}^{2}\sum\limits_{\nu=j+1}^{k}\Vert p_{\nu+1}\Vert + \lambda^{k}M_{l}h_{k}(1+l_{F}+\alpha\,h_{k})+\lambda^{k}(\alpha\,h_{k}+l_{F})\Vert \theta_{j}\Vert,
\end{aligned}
\end{equation*}
for all $ j=0,\ldots,k-1 $. Applying further Proposition~\ref{granwal1} gives us the estimate
\begin{equation*}
\begin{aligned}
\Vert p_{j+1} \Vert &\le\Big(\Vert p_{k}\Vert+\sum\limits_{i=j+1}^{k-1}(\lambda^{k}h_{k}M_{l}(1+l_{F}+\alpha h_{k})+\lambda^{k}(\alpha\,h_{k}+l_{F})\Vert \theta_{i}\Vert)\Big)\\
&\times\exp\Big(\sum\limits_{i=j+1}^{k-1}((i-j-1)\alpha h_{k}^{2}+h_{k}l_{F}+2\alpha h_{k}^{2})\Big),
\end{aligned}
\end{equation*}
valid whenever $ j=0,\ldots,k-2 $. Normalizing the discrete nontriviality condition \eqref{discr:eq1} yields
\begin{equation}\label{n:e}
\lambda^{k}+\Vert p_{k}\Vert =1,
\end{equation}
which ensures that all the sequential terms in \eqref{n:e} are uniformly bounded. Hence
\begin{equation}\label{p}
\begin{aligned}
\Vert p_{j+1} \Vert &\le\Big(1+\sum\limits_{i=j+1}^{k-1}(h_{k}M_{l}(1+l_{F}+\alpha h_{k})+(\alpha\,h_{k}+l_{F})\Vert \theta_{i}\Vert)\Big)\\
&\times\exp\Big(\sum\limits_{i=j+1}^{k-1}((i-j-1)\alpha h_{k}^{2}+h_{k}l_{F}+2\alpha h_{k}^{2})\Big).
\end{aligned}
\end{equation}
On the other hand, we have for all $ j=0,\ldots,k-2 $ that
\begin{equation}\label{p1}
\begin{aligned}
\sum\limits_{i=j+1}^{k-1}h_{k}M_{l}(1+l_{F}+\alpha\,h_{k})&=(k-j-1)h_{k}M_{l}(1+l_{F}+\alpha\,h_{k})\\
&\leq kh_{k}M_{l}(1+l_{F}+\alpha\,h_{k})=TM_{l}(1+l_{F}+\alpha\,h_{k}),
\end{aligned}
\end{equation}
\begin{equation}\label{p2}
\begin{aligned}
\sum\limits_{i=j+1}^{k-1}\Vert \theta_{i}\Vert\leq \sum\limits_{i=0}^{k-1}\Vert \theta_{i}\Vert\leq\sum\limits_{i=0}^{k-1}\int\limits_{t_{i}}^{t_{i+1}}\Big\Vert\dfrac{\bar{x}_{i+1}-\bar{x}_{i}}{h_{k}}-\dot{\bar{x}}(t)\Big\Vert\,\mathrm{d}t=\int_{0}^{T}\lVert\dot{\bar{x}}^{k}(t)-\dot{\bar{x}}(t)\rVert\,\mathrm{d}t:=\nu_{k}
\end{aligned}
\end{equation}
and thus get that $ \nu_{k}\rightarrow 0 $ as $ k\rightarrow\infty $. Furthermore, we have
\begin{equation}\label{p3}
\begin{aligned}
\sum\limits_{i=j+1}^{k-1}&\big((i-j-1)\alpha h_{k}^{2}+h_{k}l_{F}+2\alpha h_{k}^{2}\big)=\dfrac{k-j-1}{2}(k-j-2)\alpha h_{k}^{2}+(k-j-1)(h_{k}l_{F}+2\alpha h_{k}^{2})\\
&=\dfrac{(k-j-1)^{2}-(k-j-1)}{2}\alpha h_{k}^{2}+(k-j-1)(h_{k}l_{F}+\alpha h_{k}^{2})\\
&\leq \dfrac{(k-j-1)^{2}}{2}\alpha h_{k}^{2}+(k-j-1)(h_{k}l_{F}+2\alpha h_{k}^{2})\\
&\leq k^{2}\alpha h_{k}^{2}+k(h_{k}l_{F}+2\alpha h_{k}^{2})=\alpha T^{2}+Tl_{F}+2\alpha T^{2}=T(3\alpha T+l_{F}),
\end{aligned}
\end{equation}
which allows us to deduce from \eqref{p}-\eqref{p3} the desired estimates
\begin{equation}\label{es:p}
\begin{aligned}
\big\Vert p_{j+1}\big\Vert\le\big(1+TM_{l}(1+l_{F}+\alpha\,h_{k}\big)+\big(\alpha\,h_{k}+l_{F})\nu_{k}\big)\exp(T(3\alpha T+l_{F}))\leq C
\end{aligned}
\end{equation}
that hold for all $ j=0,\ldots,k-2 $ and large $k$ with some constant $ C>0 $.\\[1ex]
\textbf{Step~2.} {\em Construction and convergence of adjoint arcs.}  For each $k\in\N$, construct first the piecewise-constant functions $\th^k$ on $ [0,T] $ by
\begin{equation}\label{thetak}
\theta^{k}(t):=\theta^{k}_{j}/h_{k}\quad\text{for}\quad t\in[t^{k}_{j},t^{k}_{j+1}),\;\;j=0,\ldots,k-1,
\end{equation}
where $ \theta^{k}_{j}$ are defined in \eqref{theta}. It follows from these constructions and the convergence result of
Theorem~\ref{disc:conver} that the sequence $\{\theta^k(\cdot)\}$ converges to $ 0 $ strongly in $L^{2}([0,T],\mathbb{R}^{n})$, and hence we have that $\theta^{k}(t)\to 0$ as $k\to\infty$ for a.e. $t\in[0,T)$ along a subsequence (without relabeling). By \eqref{n:e}, suppose without loss of generality that there exists $\lm\ge 0$ with $\lm^k\to\lm$ as $k\to\infty$. Taking further $p^k_j$ with $ j=0,\ldots,k $ from Theorem~\ref{discr:N.C}, construct the adjoint functions $ p^{k}(\cdot) $ on $ [0,T] $ as the piecewise-linear extensions of $ p^{k}_{j} $. Then \eqref{es:p} means that $ p^{k}(\cdot) $ are uniformly bounded on $ [0,T] $. Moreover, using \eqref{e:p_prime}, \eqref{e:wv}, \eqref{n:e}, and \eqref{es:p} ensures that for all $ t\in(t_{j},t_{j+1}) $ and $ j=0,\ldots,k-1 $ we get
\begin{equation*}
\begin{aligned}
\Vert \dot{p}^{k}(t)\Vert &\leq M_{l}(l_{F}+1)+(2h_{k}^{-1}\Vert \mu_{j} \Vert +l_{F})\Vert p_{j+1}\Vert+h_{k}^{-1}\Vert\mu_{j}\Vert(M_{l}+\Vert\theta^{k}(t)\Vert)\\
&+h_{k}^{-1}\sum\limits_{\nu=j+1}^{k}\Vert\tilde{\xi}^{j}_{\nu}\Vert\cdot\Vert p_{\nu+1}\Vert +l_{F}\Vert\theta^{k}(t)\Vert\\
&\leq M_{l}(l_{F}+1)+(2\alpha h_{k} +l_{F})\Vert p_{j+1}\Vert+\alpha h_{k}(M_{l}+\Vert\theta^{k}(t)\Vert) +\alpha h_{k}\sum\limits_{\nu=j+1}^{k}\Vert p_{\nu+1}\Vert +l_{F}\Vert\theta^{k}(t)\Vert\\
&\leq M_{l}(l_{F}+1)+(2\alpha h_{k} +l_{F})\max\{C,1\}+\alpha h_{k}(M_{l}+\Vert\theta^{k}(t)\Vert)+\alpha T\max\{C,1\}+l_{F}\Vert\theta^{k}(t)\Vert\\
&=M_{l}(l_{F}+1)+(2\alpha h_{k} +l_{F}+\alpha T)\max\{C,1\}+\alpha h_{k}(M_{l}+\Vert\theta^{k}(t)\Vert)+l_{F}\Vert\theta^{k}(t)\Vert,
\end{aligned}
\end{equation*}
which brings us to the  estimate
\begin{equation*}
\begin{aligned}
\int_{0}^{T}\Vert \dot{p}^{k}(t)\Vert^{2}\,\mathrm{dt}&\leq 3TM^{2}_{l}(l_{F}+1+\alpha h_{k})^{2}+3T(\alpha h_{k} +l_{F}+\alpha T)^{2}(\max\{C,1\})^{2}\\
&+3(\alpha h_{k}+l_{F})^{2}\int_{0}^{T}\Vert\theta^{k}(t)\Vert^{2}\,\mathrm{dt}.
\end{aligned}
\end{equation*}
Employing this estimate and the strong convergence of $\{\theta^k(\cdot)\}$  to $ 0 $ in $L^{2}([0,T],\mathbb{R}^{n})$ implies that the sequence $ \{\dot{p}^{k}(\cdot)\} $ is bounded in $L^{2}([0,T],\mathbb{R}^{n})$, which yields the boundedness of the sequence $ \{p^{k}(\cdot)\} $ in $ W^{1,2}([0,T],\mathbb{R}^{n}) $ and hence its weak compactness in this space. Then there exists some $p(\cdot)\in
W^{1,2}([0,T],\mathbb{R}^{n})$ such that $p^{k}(\cdot)\to p(\cdot)$ uniformly in $ [0,T] $ and $\dot{p}^{k}(\cdot)\to \dot{p}(\cdot)$ weakly in  $ L^{2}([0,T],\mathbb{R}^{n}) $ along a subsequence (without relabeling).\\[1ex]
\textbf{Step~3}. {\em Construction and convergence of dual Volterra functions.} For each $k\in\N$, define the piecewise-constant functions $\sigma^{k}(\cdot)$, $ \psi^{k}(\cdot) $, and $ \gamma^{k}(\cdot)$ by
\begin{equation*}
\sigma^{k}(\tau):=h_{k}^{-1}\sum\limits_{\nu=j+1}^{k}\tilde{\xi}^{j}_{\nu}p_{\nu+1}\;\;\mbox{ if }\;\;\tau\in[t^{k}_{j},t^{k}_{j+1}),\;j=0,\ldots,k-1,\;\mbox{ with }\;\sigma^{k}(t^{k}_{k}):=0,
\end{equation*}
\begin{equation*}
\psi^{k}(\tau):=2h_{k}^{-1}\mu_{j}p_{j+1}\;\mbox{ if }\;\;\tau\in[t^{k}_{j},t^{k}_{j+1}),\;j=0,\ldots,k-1,\;\mbox{ with }\;\psi^{k}(t^{k}_{k}):=0,
\end{equation*}
\begin{equation*}
\gamma^{k}(\tau):=h_{k}^{-1}\mu_{j}\lambda^{k}(v_{j}+h_{k}^{-1}\theta_{j})\;\mbox{ if }\;\;\tau\in[t^{k}_{j},t^{k}_{j+1}),\;j=0,\ldots,k-1,\;\mbox{ with }\;\gamma^{k}(t^{k}_{k}):=0,
\end{equation*}
based on those appeared in the Volterra-type discrete optimality conditions of Theorem~\ref{discr:N.C}, and then also consider the auxiliary functions
\begin{equation}\label{vtheta}
\vartheta^{k}(t):=\max\big\{t^{k}_{j}\:\big|\;t^{k}_{j}\leq t,\,\,0\leq j\leq k\big\}\;\mbox{ for all }\;t\in[0,T],\,\,\,k\in \mathbb{N}.
\end{equation}
For $ \tau\in[t^{k}_{j},t^{k}_{j+1}) $, $ j=0,\ldots,k-1 $, the assumptions in $ (\mathcal{H}^{g}_{2}) $ and  estimate \eqref{es:p} imply that
\begin{equation*}
\begin{aligned}
\Vert\psi^{k}(\tau)\Vert\leq \dfrac{2C}{h_{k}}\int_{t_{j}}^{t_{j+1}}\Big\{\int_{t_{j}}^{t}\big\Vert\nabla\,g\big(t,s,\bar{x}^{k}(\vartheta^{k}(s))\big)^*\big\Vert\,\mathrm{d}s\Big\}\mathrm{d}t\leq (\alpha C h_{k})/2\rightarrow 0\;\;\text{as}\;\;k\rightarrow\infty,
\end{aligned}
\end{equation*}
ensuring the convergence $\psi^{k}(\tau)\to 0$ as $k\to\infty$ for all $ \tau\in[0,T]$. It follows from \eqref{e:wv} and \eqref{thetak} that
\begin{equation*}
\begin{aligned}
\int_{0}^{T}\Vert\gamma^{k}(\tau)\Vert^{2}\,\mathrm{d}\tau\leq \alpha^{2} h_{k}^{2}(TM_{l}^{2}+\int_{0}^{T}\Vert \theta^{k}(\tau)\Vert^{2}\,\mathrm{d}\tau)\rightarrow 0\;\;\text{as}\;\;k\rightarrow\infty.
\end{aligned}
\end{equation*}
Therefore, the sequence $\{\gamma^k(\cdot)\}$ converges to $ 0 $ strongly in $L^{2}([0,T],\mathbb{R}^{n})$, and hence we have that $\gamma^{k}(\tau)\to 0$ as $k\to\infty$ for a.e. $\tau\in[0,T]$ along a subsequence (without relabeling). Employing \eqref{vtheta} and the construction of $ \tilde{\xi}^{j}_{i} $ as given in \eqref{e:xii}, we have
\begin{equation}\label{e:sigma}
\begin{aligned}
\sigma^{k}(\tau)=\dfrac{1}{h_{k}}\sum\limits_{\nu=j+1}^{k-1}\int\limits_{t_{\nu}}^{t_{\nu+1}}\Big\{\int_{t_{j}}^{t_{j+1}}\nabla\,g\big(t,s,\bar{x}^{k}(\vartheta^{k}(s))\big)^*\,\mathrm{d}s\Big\}p^{k}\left(\vartheta^{k}(t)+h_{k}\right)\mathrm{d}t,
\end{aligned}
\end{equation}
for all $ \tau\in[t^{k}_{j},t^{k}_{j+1}) $, $ j=0,\ldots,k-1$. Consider further the mappings
\begin{equation*}
\tilde{g}:[0,T]\longrightarrow L^{2}([0,T],\mathbb{R}^{n})\;\mbox{ with }\; t\longmapsto \tilde{g}_{t}(\cdot):=\nabla\,g\big(t,\cdot,\bar{x}(\cdot)\big)^*,
\end{equation*}
\begin{equation*}
\tilde{g}^{k}:[0,T]\longrightarrow L^{2}([0,T],\mathbb{R}^{n})\;\mbox{ with }\;
t\longmapsto \tilde{g}^{k}_{t}(\cdot):=\nabla\,g\big(t,\cdot,\bar{x}^{k}(\vartheta^{k}(\cdot))\big)^*,
\end{equation*}
\begin{equation*}
\mu_{k} : L^{2}([0,T],\mathbb{R}^{n})\longrightarrow L^{2}([0,T],\mathbb{R}^{n})\;\mbox{ with }\;
y\longmapsto \mu_{k}\,y:=\dfrac{1}{t^{k}_{j+1}-t^{k}_{j}}\int\limits_{t^{k}_{j}}^{t^{k}_{j+1}}y(s)\,\mathrm{d}s.
\end{equation*}
Fixing $ t\in [0,T] $ and using the notation $ L^{2}:=L^{2}([0,T],\mathbb{R}^{n}) $, we get for all $ \tau\in[0,T]$ that
\begin{equation*}
\begin{aligned}
&\Vert (\mu_{k}\,\tilde{g}^{k}_{t})(\tau)-\tilde{g}_{t}(\tau) \Vert=\Big\Vert\dfrac{1}{h_{k}}\int_{t_{j}}^{t_{j+1}}\nabla\,g\big(t,s,\bar{x}^{k}(\vartheta^{k}(s))\big)^*\,\mathrm{d}s-\nabla\,g\big(t,\tau,\bar{x}(\tau)\big)^*\Big\Vert \\
&\le\dfrac{1}{h_{k}}\int_{t_{j}}^{t_{j+1}}\big\Vert\nabla\,g(t,s,\bar{x}^{k}(\vartheta^{k}(s))\big)^*-\nabla\,g\big(t,s,\bar{x}(s)\big)^*\big\Vert\,\mathrm{d}s +\Big\Vert \dfrac{1}{h_{k}}\int_{t_{j}}^{t_{j+1}}\nabla\,g\big(t,s,\bar{x}(s)\big)^*\,\mathrm{d}s-\nabla\,g\big(t,\tau,\bar{x}(\tau)\big)^*\Big\Vert\\
&\le\sup\limits_{s\in I}\big\Vert\nabla\,g\big(t,s,\bar{x}^{k}(\vartheta^{k}(s))\big)^*-\nabla\,g\big(t,s,\bar{x}(s)\big)^*)\big\Vert + \Big\Vert \dfrac{1}{h_{k}}\int_{t_{j}}^{t_{j+1}}\nabla\,g\big(t,s,\bar{x}(s)\big)^*\,\mathrm{d}s-\nabla\,g\big(t,\tau,\bar{x}(\tau)\big)^{*}\Big\Vert,
\end{aligned}
\end{equation*}
and therefore arrive at the estimate
\begin{equation*}
\Vert \mu_{k}\,\tilde{g}^{k}_{t}-\tilde{g}_{t} \Vert^{2}_{L^{2}}\leq 2T\Big(\sup\limits_{s\in I}\big\Vert\nabla\,g\big(t,s,\bar{x}^{k}(\vartheta^{k}(s))\big)^*-\nabla\,g\big(t,s,\bar{x}(s)\big)^*\big\Vert\Big)^{2}+2\Vert \mu_{k}\,\tilde{g}_{t}-\tilde{g}_{t} \Vert_{L^{2}}^{2}.
\end{equation*}
It follows from the uniform convergence $\bar{x}^{k}(\cdot)\to\bar{x}(\cdot)$ on $[0,T]$ as $k\to\infty$ and the continuous differentiability of $g(t,\tau,\cdot)$ with the usage of Lemma~\ref{mu:lemma} that
\begin{equation}\label{con:vars}
\mu_{k}\,\tilde{g}^{k}_{t}:=\dfrac{1}{h_{k}}\int_{t_{j}}^{t_{j+1}}\nabla\,g\big(t,\tau,\bar{x}^{k}(\vartheta^{k}(\tau))\big)^*\,\mathrm{d}\tau\longrightarrow \tilde{g}_{t}(\cdot):=\nabla\,g\big(t,\cdot,\bar{x}(\cdot)\big)^*\;\;\text{in}\;\;L^{2}([0,T],\mathbb{R}^{n}),
\end{equation}
for any $ t\in [0,T] $. Fixing any $ t\in [t_{\nu},t_{\nu+1}) $ with $ \nu=0,\ldots,k-1 $, we have for all $ \tau\in[0,T] $ that
\begin{equation*}
\begin{aligned}
\big\Vert(\mu_{k}\,\tilde{g}^{k}_{t})(\tau)\,p^{k}\big(\vartheta^k(t)+h_{k}\big)-\tilde{g}_{t}(\tau)\,p(t)\big\Vert &\le\big\Vert (\mu_{k}\,\tilde{g}^{k}_{t})(\tau)-\tilde{g}_{t}(\tau)\big\Vert\cdot\big\Vert p^{k}\big(\vartheta^k(t)+h_{k}\big)\big\Vert\\
&+\big\Vert p^{k}\big(\vartheta^k(t)+h_{k}\big)-p(t)\big\Vert\cdot\big\Vert \tilde{g}_{t}(\tau)\big\Vert,
\end{aligned}
\end{equation*}
which brings us to the integral inequality
\begin{equation*}
\begin{aligned}
\int_{0}^{T}\big\Vert (\mu_{k}\,\tilde{g}^{k}_{t})(\tau)\,p^{k}\big(\vartheta^k(t)+h_{k}\big)-\tilde{g}_{t}(\tau)\,p(t)\big\Vert^{2}\,\mathrm{d}\tau &\leq 2\Vert p^{k}\big(\vartheta^k(t)+h_{k}\big)\big\Vert^{2} \int_{0}^{T}\big\Vert (\mu_{k}\,\tilde{g}^{k}_{t})(\tau)-\tilde{g}_{t}(\tau)\big\Vert^{2}\,\mathrm{d}\tau\\
&+2\big\Vert p^{k}\big(\vartheta^k(t)+h_{k}\big)-p(t)\big\Vert^{2}\int_{0}^{T}\big\Vert \tilde{g}_{t}(\tau)\big\Vert^{2}\,\mathrm{d}\tau.
\end{aligned}
\end{equation*}
It follows from the uniform convergence $p^{k}(\cdot)\to p(\cdot)$ on $[0,T]$ and the convergence in \eqref{con:vars} as $k\to\infty$ with the usage of \eqref{es:p} and the imposed assumptions in $ (\mathcal{H}^{g}_{2})$ that
\begin{equation}\label{con:y}
y^{k}(t):=(\mu_{k}\,\tilde{g}^{k}_{t})p^{k}\left(\vartheta^k(t)+h_{k}\right)\longrightarrow \tilde{g}_{t}\,p(t)\;\;\text{in}\;\;L^{2}([0,T],\mathbb{R}^{n})
\end{equation}
as $ k\rightarrow\infty $ whenever $ t\in [t_{\nu},t_{\nu+1})$ for $ \nu=0,\ldots,k-1$. Observe that for $ \tau\in[t^{k}_{j},t^{k}_{j+1}) $ with $ j=0,\ldots,k-1 $, we can rewrite $\sigma^{k}(\tau)$ in \eqref{e:sigma} as
\begin{equation*}
\sigma^{k}(\tau)=\sum\limits_{\nu=j+1}^{k-1}\int\limits_{t_{\nu}}^{t_{\nu+1}}y^{k}(t)\,\mathrm{d}t=\int\limits_{t_{j+1}}^{T}y^{k}(t)\,\mathrm{d}t.
\end{equation*}
Employing the above constructions leads us to the relationships
\begin{equation*}
\begin{aligned}
&\Big\Vert \sigma^{k}(\tau)-\int_{\tau}^{T}\nabla\,g\big(t,\tau,\bar{x}(\tau)\big)^*p(t)\,\mathrm{d}t\Big\Vert\\
&=\Big\Vert \int\limits_{t_{j+1}}^{T}y^{k}(t)\,\mathrm{d}t- \int_{t_{j+1}}^{T}\nabla\,g\big(t,\tau,\bar{x}(\tau)\big)^*p(t)\,\mathrm{d}t\\
&-\int_{\tau}^{t_{j+1}}\nabla\,g\big(t,\tau,\bar{x}(\tau)\big)^*p(t)\,\mathrm{d}t\Big\Vert\\
&\le\int\limits_{t_{j+1}}^{T}\big\Vert y^{k}(t)-\nabla\,g\big(t,\tau,\bar{x}(\tau)\big)^*p(t)\big\Vert\,\mathrm{d}t+\int_{\tau}^{t_{j+1}}\Vert \nabla\,g\big(t,\tau,\bar{x}(\tau)\big)^*p(t)\Vert\,\mathrm{d}t\\
&\le(T-t_{j+1})^{\frac{1}{2}}\Big(\int\limits_{t_{j+1}}^{T}\big\Vert y^{k}(t)-\nabla\,g\big(t,\tau,\bar{x}(\tau)\big)^*p(t)\big\Vert^{2}\,\mathrm{d}t\Big)^{\frac{1}{2}}\\
&+(t_{j+1}-\tau)^{\frac{1}{2}}\disp\Big(\int_{\tau}^{t_{j+1}}\big\Vert\nabla\,g\big(t,\tau,\bar{x}(\tau)\big)^*p(t)\big\Vert^{2}\,\mathrm{d}t\Big)^{\frac{1}{2}},
\end{aligned}
\end{equation*}
which yield in turn the fulfillment of the estimates
\begin{equation*}
\begin{aligned}
&\int_{0}^{T}\Big\Vert\sigma^{k}(\tau)-\int_{\tau}^{T}\nabla\,g\big(t,\tau,\bar{x}(\tau)\big)^*p(t)\,\mathrm{d}t\Big\Vert^{2}\,\mathrm{d}\tau\\
&\le 2(T-t_{j+1})\int_{0}^{T}\Big\{\int\limits_{t_{j+1}}^{T}\big\Vert y^{k}(t)-\nabla\,g\big(t,\tau,\bar{x}(\tau)\big)^*p(t)\big\Vert^{2}\,\mathrm{d}t\Big\}\,\mathrm{d}\tau\\
&+2\int_{0}^{T}\Big\{(t_{j+1}-\tau)\int_{\tau}^{t_{j+1}}\big\Vert\nabla\,g\big(t,\tau,\bar{x}(\tau)\big)^*p(t)\big\Vert^{2}\,
\mathrm{d}t\Big\}\,\mathrm{d}\tau\\
&\le 2T\int_{0}^{T}\Big\{\int\limits_{0}^{T}\big\Vert y^{k}(t)-\nabla\,g\big(t,\tau,\bar{x}(\tau)\big)^*p(t)\big\Vert^{2}\,\mathrm{d}\tau\Big\}\,\mathrm{d}t\\
&+2h_{k}\int_{0}^{T}\Big\{\int_{\tau}^{T}\big\Vert \nabla\,g\big(t,\tau,\bar{x}(\tau)\big)^*p(t)\big\Vert^{2}\,\mathrm{d}t\Big\}\,
\mathrm{d}\tau.
\end{aligned}
\end{equation*}
Tt follows from the dominated convergence theorem and the $ L^{2}-$strong convergence \eqref{con:y} that
\begin{equation*}
\int_{0}^{T}\big\Vert \sigma^{k}(\tau)-\int_{\tau}^{T}\nabla\,g\big(t,\tau,\bar{x}(\tau)\big)^*p(t)\,\mathrm{d}t\big\Vert^{2}\,\mathrm{d}\tau\longrightarrow 0\;\;\text{as}\;\;k\rightarrow\infty.
\end{equation*}
Then there exists a subsequence of $ \{\sigma^{k}(\cdot)\} $ (again, we do not relabel) such that
\begin{equation} \label{e:con:sigma}
\sigma^{k}(\tau)\rightarrow \displaystyle\int_{\tau}^{T}\nabla\,g\big(t,\tau,\bar{x}(\tau)\big)^*p(t)\,\mathrm{d}t \;\;\text{as}\;\;k\rightarrow\infty\;\mbox{ for a.e. }\;\tau\in[0,T],
\end{equation}
which is planned to be achieved at this step.\\[1ex]
\textbf{Step~4}. {\em Verifying the necessary optimality conditions.} We start with the Volterra condition \eqref{vol-cone}, which is equivalent to \eqref{vol-cod} due to the coderivative definition \eqref{cod}. For all $ \tau\in(t^{k}_{j},t^{k}_{j+1}) $ as $ j=0,\ldots,k-1 $, rewrite \eqref{discr:eq3} as follows{
\begin{equation}\label{e:ed}
\begin{aligned}
&\dot{p}^{k}(\tau)+\psi^{k}(\tau)-\gamma^{k}(\tau)+\sigma^{k}(\tau)\in\big\{u\;\big|\;\big(u,p^{k}(\vartheta^{k}(\tau)+h_{k})-\lambda^{k}\theta^{k}(\tau)\big)\\
&\in\lambda^{k}\partial\,l\big(\vartheta^{k}(\tau),\bar{x}^{k}(\vartheta^{k}(\tau)),\dot{\bar{x}}^{k}(\tau)\big)+ N\big((\bar{x}^{k}(\vartheta^{k}(\tau)),\dot{\bar{x}}^{k}(\tau)-\bar{w}^{k}(\tau));\mathrm{gph}\,F(\vartheta^{k}(\tau),\cdot)\big)\big\},
\end{aligned}
\end{equation}
where the piecewise-constant functions $ \bar{w}^{k}: [0,T]\longrightarrow\mathbb{R}^{n} $ is given by
\begin{equation*}
\bar{w}^{k}(0):=0\quad\text{and}\quad\bar{w}^{k}(\tau):=\bar{w}^{k}_{j}\;\mbox{ if }\;\tau\in(t^{k}_{j},t^{k}_{j+1}],\;\;j=0,\ldots,k-1,
\end{equation*}
and the vectors $ \bar{w}^{k}_{j} $ are defined in \eqref{eqq:wkj}. Arguing similarly to the proof of Theorem~\ref{strong approximation} gives us the strong convergence $ \bar{w}^{k}(\cdot)\longrightarrow \bar{w}(\cdot)$ where $ \bar{w}(\tau)=\int\limits_{0}^{\tau}g(\tau,t,\bar{x}(t))\,\mathrm{d}t  $ in the $L^{2}$-norm topology as $ k\to\infty $. Then there exists a subsequence of $ \{\bar{w}^{k}(\cdot)\} $ (no relabeling) such that $ \bar{w}^{k}(\tau)\longrightarrow \int\limits_{0}^{\tau}g(\tau,t,\bar{x}(t))\,\mathrm{d}t $ almost everywhere on $[0,T] $. Mazur's theorem provides a convex combination of $ \dot{p}^{k}(\tau) $ that converges to $ \dot{p}(\tau) $ for a.e. $ \tau\in[0,T] $. Passing to the limit in \eqref{e:ed}
and employing the convergence of $ \theta^{k}(\tau)\longrightarrow 0 $, $ \dot{\bar{x}}^{k}(\tau)\longrightarrow \dot{\bar{x}}(\tau) $, and $ \vartheta^{k}(\tau)\longrightarrow \tau $ for a.e., $ \tau\in[0,T]$, with the usage of \eqref{e:con:sigma} as well as $ (\mathcal{A}^{l}) $ and $ (\mathcal{A}^{F}) $, leads us to the Volterra inclusion \eqref{vol-cone}. Taking the limit in \eqref{n:e} as $k\to\infty$ yields the nontriviality condition \eqref{nontriv}.

It remains to verify the transversality condition \eqref{tran}. To furnish this, observe first that
\begin{equation}\label{con:varphi}
\lambda^{k}\partial\varphi(\bar{x}^{k}_{k})\rightarrow \lambda\,\partial\varphi\big(\bar{x}(T)\big)\quad\text{as}\;\;k\rightarrow\infty.
\end{equation}
by the subdifferential robustness. Furthermore, the set $ \Omega_{k} $ in \eqref{e:5.1} is clearly represented as
\begin{equation*}
\Omega_{k}=\big\{x_{k}\in\mathbb{R}^{n}\;\big|\;\mathrm{dist}\,(x_{k};\Omega)\leq \zeta_{k}\big\},
\end{equation*}
via the distance function. Passing to the limit in \eqref{discr:eq2} as $ k\rightarrow\infty $ and using the normal cone representation in \eqref{cone:dist} together with \eqref{con:varphi}, we justify \eqref{tran} and thus complete the proof. $\h$\vspace*{-0.2in}

\section{Conclusions and Future Research}\label{sec:concl}\vspace*{-0.1in}

The paper develops the method of discrete approximations to derive necessary optimality conditions for a new class of dynamical systems governed by integro-differential inclusion. The developed machinery of discrete approximations is of its own theoretical and numerical interest, while it helps us to establish a novel type of optimality conditions that are coined here as ``generalized Volterra."

In our future research, we plan to further develop the method of discrete approximations, in its both theoretical and numerical aspects, to cover various classes of controlled sweeping processes governed by integro-differential inclusions and to employ the obtained results in deriving novel necessary optimality conditions of the generalized Volterra type. Moreover, our goals will be applying the new conditions to optimization of practical models that can be descried by controlled, both Lipschitzian and discontinuous, integro-differential inclusions that particularly appear in electronic modeling.
\end{document}